%% file: ex_article.tex
\documentclass[hidelinks,onefignum,onetabnum]{siamart220329}

\usepackage{amsfonts}
\input{ex_shared}

\ifpdf
\hypersetup{
  pdftitle={Solving McKean-Vlasov Equation by deep learning particle method},
  pdfauthor={Jingyuan Li ,Wei Liu}
}
\fi

\begin{document}

\maketitle

\begin{abstract}
We introduce a novel meshless simulation method for the McKean-Vlasov Stochastic Differential Equation (MV-SDE) utilizing deep learning, applicable to both self-interaction and interaction scenarios. Traditionally, numerical methods for this equation rely on the interacting particle method combined with techniques based on the It\^o-Taylor expansion. The convergence rate of this approach is determined by two parameters: the number of particles $N$ and the time step size $h$ for each Euler iteration. However, for extended time horizons or equations with larger Lipschitz coefficients, this method is often limited, as it requires a significant increase in Euler iterations to achieve the desired precision $\epsilon$. To overcome the challenges posed by the difficulty of parallelizing the simulation of continuous interacting particle systems, which involve solving high-dimensional coupled SDEs, we propose a meshless MV-SDE solver grounded in Physics-Informed Neural Networks (PINNs) that does not rely on the propagation of chaos result. Our method constructs a pseudo MV-SDE using It\^o calculus, then quantifies the discrepancy between this equation and the original MV-SDE, with the error minimized through a loss function. This loss is controlled via an optimization algorithm, independent of the time step size, and we provide an error estimate for the loss function. The advantages of our approach are demonstrated through corresponding simulations.
\end{abstract}

\begin{keywords}
McKean-Vlasov equation, Physics-Informed Neural Networks, Particle method, Meshless solver
\end{keywords}


\section{Introduction}
In recent years, McKean-Vlasov Stochastic Differential Equations (MV-SDEs) have garnered significant attention across various disciplines, including finance, biology, and chemistry \cite{Baladron2012MeanfieldDA,lasry2018mean,carmona2018probabilistic}.MV-SDEs are crucial as macroscopic limit equations for mean-field interacting particle systems, making them foundational in fields such as statistical physics and machine learning. Despite their widespread applicability, two primary challenges hinder their effective use: efficient computation of the stationary distribution and the development of accurate, efficient numerical solvers. While the former can be addressed by allowing the simulation time to approach infinity, the latter remains a persistent issue, particularly for long time horizons.

Numerical solvers for MV-SDEs typically rely on interacting particle systems or self-interacting diffusions, which necessitate solving high-dimensional coupled equations. These solvers aim to extract distributional information, such as moments or density functions. Recent advancements include Monte Carlo methods \cite{rached2022single,Rached2022MultilevelIS,veretennikov2006ergodic} and self-interaction techniques \cite{du2023empirical} that utilize sampling to estimate target distributions. A classical approach involves approximating particle systems using Euler-Maruyama method \cite{Kloeden1977TheNS,reisinger2023convergence,Bogachev2019OnCT}. In the context of interacting particle systems, the distributions of MV-SDE solutions are often approximated by empirical measures of the particles, resulting in high-dimensional coupled SDE systems. The convergence rate of this approach depends on two key factors \cite{Antonelli2002RateOC}: the number of particles $N$, which controls the deviation between the empirical and true distributions typically through the propagation of chaos \cite{Salinas2020ATF}, and the time step size
$h$, determined by the iteration scheme, such as the Euler method, truncated Euler method, or other variants \cite{Bao2021FirstorderCO,He2022AnEE,Liu2022ParticleMA}. Although existing methods provide theoretical bounds, they often suffer from computational inefficiencies due to the complexity of solving these high-dimensional systems.

Another class of numerical solvers employs self-interacting diffusions, where the true distribution is approximated by the occupied measure of the diffusion, analogous to the empirical measure in particle systems. While self-interacting diffusion methods show promise for simulating single particles, they are limited by the complexity of numerical solvers and the need to store a large number of historical samples. Beyond these two approaches, another method involves solving the nonlinear Fokker-Planck equation \cite{Barbu2018FromNF} to gain distributional insights. However, this approach does not provide individual particle trajectory and can be computationally demanding.

Given the increasing integration of deep learning in scientific computing, we propose a novel deep learning-based solver for MV-SDEs to address these challenges. Our primary contribution is the direct approximation of particle systems within a continuous framework using Deep Neural Networks (DNNs). We construct a pseudo MV-SDE using It\^o calculus, with differential terms in the drift and diffusion coefficients computed via automatic differentiation. The loss function is defined as the $L^2$ error between the pseudo MV-SDE and the original MV-SDE, which is minimized using an optimization algorithm. This approach is applied to approximate interacting particle systems, aiming to enhance the accuracy of numerical solvers. Similarly, we utilize this method to approximate self-interacting diffusions, thereby obtaining stationary distributions without the need to store historical path samples.

Our method offers an alternative and improvement over traditional numerical solvers, where the iteration format is often predetermined, and typically only a single time point of the Brownian motion is simulated to compute the particle's position in the next step. In contrast, our approach requires simulating the entire trajectory of the Brownian motion over the time interval in advance. The primary advantages of our method include its meshless nature and its capability for GPU acceleration. Furthermore, the method is versatile, accommodating various types of noise, including fractional Brownian motion, by simply adjusting the corresponding It\^o formula.

We provide an error estimate for our approximation and validate our method through a series of experiments. For interacting particle systems, our approach achieves higher accuracy compared to traditional solvers with fixed step sizes while maintaining a meshless framework. The effectiveness of our method is further demonstrated under fractional Brownian motion noise. For self-interacting diffusions, our approach successfully approximates the stationary distribution along a single trajectory.

The rest of this paper is  organized as follows. In Section 2, we provide the backgrounds on MV-SDEs and introduce our method along with its fundamental approach. We offer detailed explanations and analyses and specifically deal with two distinct scenarios: interacting particle approximation and self-interacting cases. In Section 3, we perform detailed error analysis for these scenarios, including the error evaluations in the sense of $L^2$ and Wasserstein distance. Section 4 is devoted to presenting our algorithm in detail and describing a series of experiments conducted to validate the effectiveness of our approach.

\section{Introduction of the methodology}

In this section, we introduce the framework and propose our method. The goal of this paper is to propose a numerical solver for the following MV-SDE:
\begin{equation}\label{eq1}
\begin{aligned}
\left\{\begin{array}{l}
d X_{t}=b\left(t, X_{t}, \mu_{t}\right) d t+\sigma\left(t, X_{t}, \mu_{t}\right) d W_{t} \\
\mu_{t} = \mathcal{L}(X_{t}) \text{ is the distribution of $X_t$}
\\
\mathcal{L}(X_0)=\mathcal{L}_0 \text{ is the distribution of $X_0$}, \quad\quad t \in [0, T]
\end{array}\right.
\end{aligned}
\end{equation}
where $X_t$ is defined in some probability space $\left(\Omega, \mathcal{F},\left(\mathcal{F}_{t}\right)_{t \geq 0}, \mathbb{P}\right)$,  $X_0$ is the initial value, and $\{W_t,t\ge 0\}$ is a standard Brownian motion. Let $\mathcal{P}_p\left(\mathbb{R}^d\right)$ denote the space of probability measures on $\mathbb{R}^d$ with finite moments of order $p$. The drift term coefficient $b$ is a Borel measurable, $\mathbb{R}^d$-valued function defined on $[0, T] \times \mathbb{R}^d \times \mathcal{P}(\mathbb{R}^d)$, and diffusion term coefficient $\sigma: [0, T] \times \mathbb{R}^{d} \times \mathcal{P}\left(\mathbb{R}^{d}\right) \to  \mathbb{M}_{d, q}(\mathbb{R})$ is a matrix-valued function with $\mathbb{M}_{d, q}(\mathbb{R})$ being the space of $d \times q$ real matrices.

The corresponding interacting particle system takes the following form: \begin{equation}\label{eq2}
d X^{n,N}_{t}=b\left(t, X^{n,N}_{t}, \mu^{(N)}_{t}\right) d t+\sigma\left(t, X^{n,N}_{t},\mu^{(N)}_{t}\right) d W_{t}^n, \quad\quad t \in [0, T], \quad\quad 1\le n\le N
\end{equation}
where $\mu^{(N)}_{t}:=\frac{1}{N} \sum_{n=1}^{N} \delta_{{X}_t^{n, N}}$ is the empirical measure of the particles ${X}_t^{n, N}, 1\le n\le N$, and $X_0^{n,N}$ are i.i.d. random variables with the same distribution as $X_0$. Due to the pairwise interactions of the particles, ${X}_t^{n, N}, 1\le n\le N$ are usually not independent but have the same distribution. Under some assumptions (such as Lipschitz continuity and linear growth of the coefficients $b$ and $\sigma$), the empirical measure $\mu^{(N)}_t$ (equivalently, the law of each particle) will converge to the distribution $\mu_t$ of the solution $X_t$ of the McKean-Vlasov equation \cref{eq1}, as the number of particles $N$ tends to infinity (so-called Propagation of Chaos). This is often done by considering the following non-interacting (Auxiliary) particle system:
\begin{equation}\label{eq3}
d \tilde{X}^{n}_{t}=b\left(t, \tilde{X}^{n}_{t}, \mu^{n}_{t}\right) d t+\sigma\left(t, \tilde{X}^{n}_{t}, \mu^{n}_{t}\right) d W_{t}^n,  \quad\quad t \in [0, T],\quad\quad 1\le n\le N
\end{equation}
where $\mu^{n}_{t}$ is the distribution of $\tilde{X}^n_t$, and $\tilde{X}_0^{n}$ are i.i.d. random variables with the same distribution as $X_0$. Indeed, $\tilde{X}^{n}_{t}, 1\le n\le N$ are independent copies of the solution $X_t$ in \cref{eq1}. Then one may use coupling methods or functional inequalities to estimate the Wasserstein distance between the distributions of independent particles $\tilde{X}^{n}_{t}$ and interacting particles $X^{n,N}_{t}$, which converges to $0$ as the number $N$ of the particles tends to infinity, see \cite{sznitman1991topics,durmus2020elementary}.

This provides a numerical solver for the McKean-Vlasov equation. If we apply the Euler-Maruyama method directly to the McKean-Vlasov equation, we will encounter the problem of unknowing the distribution $\mu_t$ of $X_t$.  However, due to propagation of chaos, we only need to apply the Euler-Maruyama method to the interacting particle system, with the empirical measure $\mu^{(N)}_t$ of interacting particles replacing $\mu_t$. The following is a Euler-Maruyama discretization of the interacting particle system (\ref{eq2}):
\begin{eqnarray}
\label{eq4}
\quad\quad\left\{\begin{array}{l}
X_{t_{m+1}}^{n, N}=X_{t_{m}}^{n, N}+h b\left(t_{m}, {X}_{t_{m}}^{n, N}, \hat{\mu}_{t_{m}}^{N}\right)+\sigma\left(t_{m}, X_{t_{m}}^{n, N}, \hat{\mu}_{t_{m}}^{N}\right)\left(W_{t_{m+1}}^{n}-W_{t_{m}}^{n}\right) \\
\hat{\mu}_{t_{m}}^{N}:=\frac{1}{N} \sum_{n=1}^{N} \delta_{{X}_{t_{m}}^{n, N}}
\\
\end{array}\right.\text{,}
\end{eqnarray}
where $(t_0,t_1,\ldots,t_m)$ is the time point with the step size $h=t_m-t_{m-1}$, $\forall m \in 1,2,\dots,M$.  When $h$ tends to zero and $N$ tends to infinity, the Euler-discretized particle in (\ref{eq4}) converges to the MV-SDE (\ref{eq1}) with the convergence rate $\mathcal{O}\left(\frac{1}{n} + h\right)$ in the sense of $L^1$ norm, as discussed in \cite{Antonelli2002RateOC}.

{\bfseries We mainly focus on the structure of numerical algorithms, particularly those designed to approximate interacting (or self-interacting) particle systems. Essentially, our method aligns with the approach outlined in equation (\ref{eq4}), but differs in that the error is not affected by the time step size.}

\subsection{Methods for approximating interacting particle system}

The basic idea of this approach is to preserve the assumptions of exchangeability as discussed in \cite{Chaintron2022PropagationOC}, ensuring that the solution $X_t^{i,N}$ remains invariant under permutations of the particle indices, thus maintaining the symmetry of this model. The particles $X_t^{i,N}$ in (\ref{eq2}) are identically distributed but not mutually independent. Once the propagation of chaos holds, the solution of (\ref{eq2}) will approximate that of (\ref{eq1}) when $N$ is very large. In this context, an alternative approach to approximate the interacting particle system (\ref{eq2}) is considered, by representing the solution as the functional of the corresponding noise:
\begin{equation}\label{eq5}
\begin{aligned}Y_t^{n,N} = F^n(t,W^n_t)\end{aligned}\text{,}
\end{equation}
where $F^n \in C^{1,2}$ is  continuously differentiable with respect to the first variable and twice continuously differentiable with respect to the second variable. Then we can choose an appropriate structure of neural network with a nonlinear activation function to approximate $F^n$ \cite{Ryck2021OnTA}. In typical cases, one can specify a particular loss function and then use optimization algorithms such as SGD to obtain a neural network approximation of the objective function. The deep learning-based method for approximating PDE is widely used, primarily aimed at minimizing the residuals of the PDE to obtain an estimate of the solution at the respective time points, typically implemented in the form of \textit{Physics-Informed Neural Networks (PINN)} \cite{Raissi2019PhysicsinformedNN,Lagaris1997ArtificialNN}, etc. In this paper, our objective is to find, for each particle $X^{n,N}$, a function $F^n(t, W^n_t) \in C^{1,2}$ that matches the drift and diffusion terms of (\ref{eq2}). Therefore, the It\^o formula can be applied to equation (\ref{eq5}):
\begin{equation*}
dY_t^{n,N} = F^n_1(t,W_t^n)dt + \frac{1}{2}F^n_{22}(t,W_t^n)dt + F^n_2(t,W_t^n)dW_t^n\text{,}
\end{equation*}
where $F^n_1$ and $F^n_2$ denote the partial derivatives of $F$ with respect to the time variable $t$ and space variable $x$, and $F^n_{22}$ denotes the second order partial derivative of $F$ with respect to the space variable $x$.

For simplicity, we will use the following notations:
\begin{eqnarray}
    b ^{Pseudo}_{n,t}:=F^n_1(t,W_t^n)+ \frac{1}{2}F^n_{22}(t,W_t^n) \text{,}\\
    \sigma ^{Pseudo}_{n,t}:=F^n_2(t,W_t^n)   \text{.}
\end{eqnarray}
We aim to achieve the following approximation through training:
\begin{eqnarray}\label{eq6}
b ^{Pseudo}_{n,t} \approx b(t,Y^{n,N}_t,\hat \nu_t{(N)})\text{,} \\
\label{eq7}\sigma ^{Pseudo}_{n,t} \approx \sigma(t,Y^{n,N}_t,\hat \nu_t{(N)}) \text{,}
\end{eqnarray}
where $\hat{\nu}_t^{(N)}:=\frac{1}{N} \sum_{n=1}^{N}\delta_{Y^{n,N}_t}$ is the empirical measure of $Y^{n,N}_t, 1\le n\le N$. This will allow us to use $Y^{n,N}_t$ to approximate $X^{n,N}_t$ for each $1\le n\le N$, and the error between  $X^{n,N}_t$ and $Y^{n,N}_t$ can be controlled by two approximation errors in (\ref{eq6}) and (\ref{eq7}) above, which will serve as the main components of our loss function in the neural network. It is very important to note that this error is independent of the time step size $h$.

\subsection{Methods for approximating independent particle system}
{\bfseries Another numerical scheme has a smaller number of parameters and reduces the amount of computation.}
A key difference between deep learning-based methods and classical numerical solvers is the ability to compute the trajectory of Brownian motion in advance. Due to this nature of deep learning, it is feasible to simulate the MV-SDE with only one particle
 The same loss function can be employed as in previous method. Optimizing this loss function results in an approximation of the equation, involving a global adjustment that highlights the meshless nature of the method. In contrast, classical numerical methods necessitate complete information from the previous time step for iteration to obtain the next time step, representing a local adjustment.

\begin{figure}[h]
\centering
\includegraphics[width=0.9\textwidth]{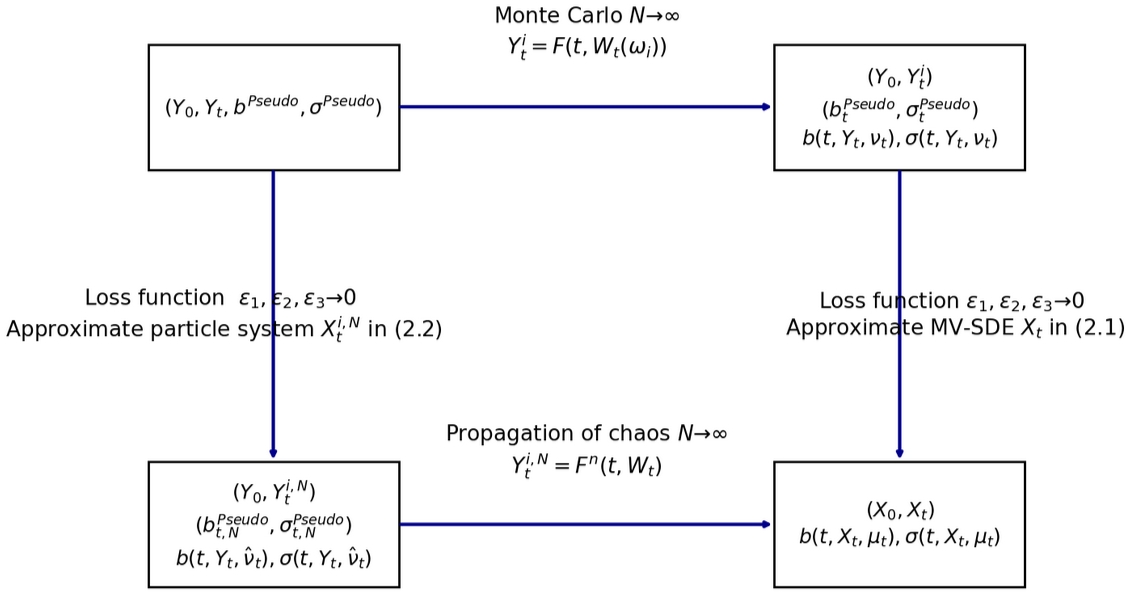}
\caption{Two different approximation methods for MV-SDE}
\label{flow}
\end{figure}

Explicitly, a neural network $F_{\theta}$ is considered to approximate the solution of equation (\ref{eq3}) directly. We can sample $N$ Brownian motion paths independently, and then get independent samples of $Y_t^n$ as follows:
\begin{equation*}Y_t^n:=Y_t(\omega^n) = F_{\theta}(t,W_t(\omega^n)) \text{ for almost every } \omega^n,   \quad 1\le n\le N.\end{equation*}
During this process, we take $N$ large enough to use the empirical distribution of $Y_t^n, 1\le n\le N$ to replace $\mu_t$. Then we choose the loss function similar to the previous method by using (\ref{eq6}) and (\ref{eq7}).
Notably, in the case of using a single $F$, since the neural network can anticipate the trajectories of different Brownian motions corresponding to the same $F$, and these $ Y^n_t$ are i.i.d., we can directly apply the results of Wasserstein distance-convergence in the i.i.d. setting \cite{Fournier2013OnTR}. \textbf{This allows us to omit the propagation of chaos, marking a theoretical breakthrough}. This is one of the features of using deep learning methods. We will discuss this in more detail in the convergence proof later.
Figure (\ref{flow}) shows two different approximation methods for MV-SDE, where the use of a single $F$ is specific to neural network methods. 

Classical methods (e.g., \ref{eq4}) first let the step size $h$ approach zero, converging to the interacting particle system, and then let $N$ tend to infinity. This process relies on the propagation of chaos. For the multi-trajectory case (i.e., with different $F$), this corresponds to an approximate particle system, consistent with the lower path in the figure, while our method for the multi-trajectory case is grid-free. In contrast, for the single-trajectory case (i.e., with the same $F$), $N$ tends to infinity first, followed by minimizing the loss, which does not rely on the propagation of chaos. This fundamental distinction highlights the unique characteristics of deep learning-based methods. However, for the single-trajectory case, the lack of consideration for particle interactions leads to a slower convergence rate when directly using Monte Carlo integration, as evidenced by the theoretical results presented later.

\subsection{Approximating stationary distributions for self-interacting systems}

Only one path is considered to approximate the distribution of the MV-SDE. Typically, self-interaction methods are employed for this purpose. This approach also facilitates the simulation of the stationary distribution of MV-SDE, specifically by employing time averaging to obtain an empirical distribution.  Under some mild conditions \cite{du2023empirical,wang2018distribution},

The stationary distributions (also referred as equilibrium states) of SDE are very important in many fields such as physics, chemistry, biology, finance, machine learning, etc. For the existence and uniqueness of stationary distributions of MV-SDEs, one can refer to \cite{guillin2022uniform,liu2021long}. If there are more than one stationary distribution, we say that the phase transition occurs.  However it is usually not to be explicitly represented or exactly computed, so we always have to do simulation and error analysis, see \cite{guillin2021kinetic}. This is also one of the reasons why the theoretical convergence rate to the stationary distribution has received much attention.

In simulation the stationary distribution, interacting particle system is often adopted to replace $\mu_t$ with the empirical distribution of the particles (spatial average). Another way to approximate $\mu_t$ is to use the empirical distribution in the sense of time average of the MVSDE itself , not the interacting particle system. The authors also estimate the convergence rate and error in \cite{du2023empirical}. Explicitly, consider
the following self-interacting process:
\begin{equation}\label{eq8}
d X_t=b\left(X_t, \kappa_{t}\right) d t+\sigma\left(X_t,\kappa_{t}\right)dW_t\text{,} \quad \kappa_{t} = \frac{1}{t}\int_{o}^{t} \delta_{X_s}ds\text{.}
\end{equation}

It is more straightforward to approximate the MV-SDE using the self-interacting diffusion (\ref{eq8}) above. This is notably useful for approximating the stationary distribution since the empirical distribution $\kappa_{t}$ will converge to the stationary distribution when $t$ tends to infinity due to the ergodicity.  In practice, delayed or discrete sampling is much common used to replace $\kappa_{t}$. Given $N\ge 1$, let $t_k=\frac{k}{N}T, k=0,1\cdot, N$. For any $t\in (0,T]$ with $t_k\le t< t_{k+1}$, let $\hat\kappa_{t}= \frac{1}{t_k} \sum_{i=1}^{k} \delta_{X_{t_i}}$, then we consider the following process
\begin{equation}\label{eq9}
d \hat{X}_t=b\left(\hat{X}_t, \hat \kappa_{t}\right) d t+\sigma\left(\hat{X}_t, \hat\kappa_{t}\right)dW_t\text{.}
\end{equation}
It has been shown in \cite{du2023empirical} that when $t$ approaches infinity, $\hat{X}_t$ will converge to the stationary distribution of equation (\ref{eq8}) under the 2-Wasserstein distance.

This type of MV-SDE solver has received little attention. A practical idea is to assume that $\hat{\kappa}_{t}$ remains constant over a period of time, and it is easily satisfied when approaching a stationary distribution. Then a discrete approximation of the equation (\ref{eq9}) can be performed using the delayed equation solver. But this method requires that the distribution dependence term is invariant with time over a period of time. This method is not clear and accordingly stores the results of a large number of samples from previous moments. In this method, this issue is mitigated since the positions of all $X_t$ can be determined in advance, facilitating the calculation of $\hat\kappa_{t}$.

Similarly as before, we consider the following direct approximation method for equation (\ref{eq9}) by applying neural network to a single particle:
\begin{equation}
Y_t = F_{\theta}(t,W_t)\text{.}
\end{equation}
The It\^o formula is used to construct the following loss function:
\begin{eqnarray}\label{eq10}
F_1(t,W_t)+ \frac{1}{2}F_{22}(t,W_t)=: b ^{Pseudo}_{n,t}  \approx b(Y_t,\hat \kappa_{t})\text{,} \\
\label{eq11}F_2(t,W_t)=: \sigma ^{Pseudo}_{t} \approx \sigma(Y_t,\hat \kappa_{t})\text{.}
\end{eqnarray}
The self-interacting Pseudo MV-SDE is constructed as follows:
\begin{equation}\label{eq12}
d Y_t=b ^{Pseudo}_{t}d t+\sigma ^{Pseudo}_{t}dW_t\text{,} \quad \hat\kappa_{t} = \frac{1}{t_k} \sum_{i=1}^{k} \delta_{Y_{t_i}}\text{.}
\end{equation}
The error estimate of this method will be provided for in next session.

\subsection{Reasons and Details}
The design of the architecture $F_{\theta}(t,w)$ is motivated by several key factors:

\textbf{Handling Uneven Data Point Distribution and Meshless Characteristics}: This method efficiently manages unevenly distributed data points by utilizing them exclusively for training the function estimation. As a result, there is no need to impose smaller time step sizes, which is typically required in traditional methods.


\textbf{Enhanced Error Control}: In conventional MV-SDE solvers, errors often grow exponentially with the time horizon $T$ and are further influenced by the Lipschitz constants of the drift and diffusion coefficients. Classical approaches require smaller time steps
$h$ to mitigate these errors, resulting in more iterations and increased computational effort, particularly over longer timescales. In contrast, our method controls error through the loss function, which is independent of the time step size. By leveraging various optimizers to minimize the loss, we provide a more efficient and flexible error control mechanism.

\textbf{Broad Applicability to Various Noise Types}: This framework can be applied to a wide range of noise types, including McKean-Vlasov SDEs driven by fractional Brownian motion \cite{He2022AnEE}. By adapting the treatment of the noise and the corresponding stochastic chain rule (It\^o formula), our method accommodates these noises and enables accurate error estimation.

\textbf{Theoretical Properties and Computational Efficiency}: From a computational perspective, when particles interact (as described in equation (\ref{eq2})), each particle requires its own neural network, optimizer, and related components. However, in non-interacting scenarios, a single neural network can be employed for all particles, which decouples their equations (as described in equation (\ref{eq3})). This leads to robust theoretical properties without the need for additional conditions beyond those of typical MV-SDE solvers. Under the propagation of chaos, our approach enhances efficiency, especially when using only one particle. Moreover, in the non-interacting case, the shared neural network significantly reduces computational overhead by eliminating the need for separate networks and optimizers for each particle. This results in substantial efficiency gains when handling large numbers of particles, bypassing the limitations imposed by propagation of chaos since all particles are independent and identically distributed (\textit{i.i.d}).

\section{Theoretical properties and error analysis}
We provide the error estimates for our method in this section .
\subsection{Notations and Definition}
 Denote $\|X\|_{p}:=\left[\mathbb{E}|X|^{p}\right]^{1 / p}\text{, } p \geq 2$ as $L^p $ norm for a random variable $X $. Let $||\sigma||$ be the Frobenius norm of the matrix $\sigma$. For any two probability measures $\mu,\nu\in\mathcal{P}_p$, define the $p$-order Wasserstein distance between $\mu$ and $\nu$ as follows:
\begin{eqnarray*}
 & W_{p}(\mu, \nu):=\left(\inf _{\pi \in \Pi(\mu, \nu)} \int_{\mathbb{R}^{d} \times \mathbb{R}^{d}}|x-y|^{p} \pi(d x, d y)\right)^{\frac{1}{p}} \\
& =\inf \left\{\left[\mathbb{E}|X-Y|^{p}\right]^{\frac{1}{p}}| X, Y:(\Omega, \mathcal{F}, \mathbb{P}) \rightarrow\left(\mathbb{R}^{d}, \mathcal{B}\left(\mathbb{R}^{d}\right)\right) \text { with } \mathcal{L}_{X}=\mu, \mathcal{L}_{Y}=\nu\right\} \text{,}
\end{eqnarray*}
where $\Pi(\mu, \nu)$  denotes the class of product measures on $\mathbb{R}^{d} \times \mathbb{R}^{d}$ with marginal $\mu$ and $\nu$.

The $L^p$ space and associated $L^p$ norm are defined as follows:
\begin{equation*}
L^p =  \left\{ f :\int |f|^pd\mu<+\infty \right\} \text{,}\|f\|_p:=(\int |f|^pd\mu)^{\frac{1}{p} }\text{.}
\end{equation*}

The Sobolev space $W^{k,p}$ is defined by
\begin{equation*}
W^{k,p} = \left\{f \in L^p :D^mf \in L^p \text{for all } |m| \leq k \right\}\text{,}
\end{equation*}
where
\begin{equation*}
D^mf = \frac{\partial^{|m|}f}{\partial^{m_1} x_1 \partial^{m_2} x_2 \dots\partial^{m_d} x_d}
\quad |m|=\sum_{i=1}^{d} m_i \text{.}
\end{equation*}

Let $\mathcal{NN}^n_m$ denote the space of fully connected neural networks with maximum width of $n$ and depth of $m$. If $F \in \mathcal{NN}^n_m$, then $F$ can be written as $F_{\theta}(t,w)$, where $\theta$ is the training parameter of the neural network. 

\subsection{Preparation Lemmas}

\begin{lemma}[The Generalized Minkowski Inequality]
\label{lemma1}
For any process $X=\left(X_{t}\right)_{t \geq 0}$, for every $p \in[2, \infty)$ and for every $T \in[0,+\infty]$,
\begin{equation*}
\left\|\int_{0}^{T} X_{t} d t\right\|_{p} \leq \int_{0}^{T}\left\|X_{t}\right\|_{p} d t .
\end{equation*}
\end{lemma}

\begin{lemma}[Burk?lder-Davis-Gundy Inequality]
\label{lemma2}
For every $p \in (0, +\infty)$, there
exists two real constants $c_p$
and $C_p>0$ such that, for every continuous local martingale $X_t$ $(t \in [0,T])$ null at 0.
\begin{equation*}
c_p\|\sqrt{\left \langle X \right \rangle_T } \|_p \leq  \left \|\sup_{t \in [0,T]}|X_t| \right\| _p\leq C_p\|\sqrt{\left \langle X \right \rangle_T }\|_p \text{,}
\end{equation*}
where \( \langle X \rangle_T \) denotes the quadratic variation of \( X_t \) up to time \( T \).
\end{lemma}

\begin{lemma}[la Gronwall Inequality]
\label{lemma3}
$f:[0,T] \to \mathbb{R^+}$ is a Borel, locally bounded, non-negative
and non-decreasing function and let $\varphi : [0,T] \to \mathbb{R^+}$  be a non-negative non-decreasing function satisfying
\begin{equation*}
f(t) \leq A\int_{0}^{t} f(s)ds +B \left(\int_{0}^{t} f^2(s)ds\right)^\frac{1}{2} +\varphi (t)  \quad \text{for any } t \in [0,T]\text{,}
\end{equation*}
where $A, B$ are two positive real constants. Then,for any $t \in [0,T]$,
\begin{equation*}f(t) \leq 2e^{(2A+B^2)t}\varphi (t)\text{.}\end{equation*}
\end{lemma}

Let $\hat\mu_{N}:=\frac{1}{N} \sum_{n=1}^{N} \delta_{X_{n}}$ be the empirical distribution of i.i.d. random variables $X_n, 1\le n\le N$ with common law $\mu$. Next we present the result on the  Wasserstein distance between $\hat\mu_{N}$ and $\mu$ \cite{Fournier2013OnTR}.
\begin{lemma}
\label{lemma4}
Let  $p > 0$  and $\mu \in \mathcal{P}_{q}\left(\mathbb{R}^{d}\right) $ for some  $q>p$.
Then there exists a positive constant $C$ only depending on $p,d,q$ such that, for all $N \geq 1$
\begin{align*}
&\mathbb{E}\left({W}_{p}^{p}\left(\hat\mu_{N}, \mu\right)\right) \leq C M_{q}^{p / q}(\mu) \times \\ &\left\{\begin{array}{ll}
N^{-1 / 2}+N^{-(q-p) / q} & \text { if } p>d / 2 \text { and } q \neq 2 p \\
N^{-1 / 2} \log (1+N)+N^{-(q-p) / q} & \text { if } p=d / 2 \text { and } q \neq 2 p \\
N^{-p / d}+N^{-(q-p) / q} & \text { if } p \in(0, d / 2) \text { and } q \neq d /(d-p)
\end{array}\right.\text{,}
\end{align*}
where  $M_{q}(\mu):=\int_{\mathbb{R}^{d}}|x|^{q} \mu(dx)$.
\end{lemma}

\subsection{Theoretical results}

\begin{hypothesis}\label{hyp1}
There exists some $p \in[2,+\infty)$ such that the following conditions (1) and (2) hold.

(1) The initial random variable $X_{0}$ satisfies $\left\|X_{0}\right\|_{p}<+\infty$.

(2) The coefficient functions $b$ and $\sigma$ are $\gamma$-H?lder continuous in $t, \gamma \in(0,1]$, and Lipschitz continuous in $x$ and in $\mu$ in the following sense: there exists a constant $L>0$ such that
\begin{eqnarray*}
\begin{aligned}
 \forall(x, \mu&) \in \mathbb{R}^{d} \times \mathcal{P}_{p}\left(\mathbb{R}^{d}\right), \forall s, t \in[0, T] \text { with } s \leq t, \\
 \quad|b(t, x, \mu)-b(s, x, \mu&)| \vee\|\sigma(t, x, \mu)-\sigma(s, x, \mu)\| \mid \leq L\left(1+|x|+W_{p}\left(\mu, \delta_{0}\right)\right)(t-s)^{\gamma},
\end{aligned}
\end{eqnarray*}
and such that
\begin{eqnarray*}
\begin{aligned}
&\forall t \in  {[0, T], \forall x, y \in \mathbb{R}^{d} \text { and } \forall \mu, \nu \in \mathcal{P}_{p}\left(\mathbb{R}^{d}\right) } \\
&|b(t, x, \mu)-b(t, y, \nu)| \vee \|\sigma(t, x, \mu)-\sigma(t, y, \nu)\| \leq L\left(|x-y|+W_{p}(\mu, \nu)\right)\text{.}
\end{aligned}
\end{eqnarray*}
Here we use the same hypothesis as \cite{Liu2022ParticleMA,liu2019optimal}.
\end{hypothesis}

Next we will present results for two cases: 1) each particle is associated with a distinct functional $F$, and the propagation of chaos is assumed to hold in this scenario; 2) all particles share the same functional $F$, and the MV-SDE is directly approximated without relying on propagation of chaos.

\begin{theorem}[Estimation of moments and Wasserstein distance bound in particle system]
Assume that Hypothesis (\ref{hyp1})  holds.
Consider a system of particles $X_t^{n,N}, t \in [0,T]$ defined in (\ref{eq2}), and let $Y_t^{n,N} = F^n(t,B_t^{n,N})$ where $F^n$ are neural networks. Denote the training time points by $t_0=0,t_1, \dots ,t_m=T \in [0,T]$. Assume that the following conditions hold for all $1\le n\le N$:
\begin{eqnarray}\label{cod1}
\|X_0^{n,N} - Y_0^{n,N} \|_p \leq \epsilon_1   \text{,}\\
\label{cod2}\int_{0}^{t} \left \|b(s,Y_s^{n,N},\hat \nu_s ) -  b ^{Pseudo}_{n,s} \right \|_pds   \leq \epsilon_2 \text{,}\\
\label{cod3}\left  \| \sqrt{\int_{0}^{t} \left  \| \sigma(s,Y_s^{n,N},\hat \nu_s) -  \sigma ^{Pseudo}_{n,s} )\right\|^2 ds}\right\|_p \leq \epsilon_3 \text{,}
\end{eqnarray}
where $ \hat \nu_t :=\frac{1}{N} \sum_{n=1}^{N}  \delta_{{Y}_{t}^{n, N}}$. Then there exist two constants $C_1$ and $C_2$ depending on $L,T,N,d,\epsilon_1,\epsilon_2,\epsilon_3$, such that
\begin{equation*}
\left \| \sup _{u \in [0,t]}|X_u^{n,N}-Y_u^{n,N}|
\right \|_p\leq
 C_1
\end{equation*}
and
\begin{equation*}
 \sup _{u  \in [0,t]}\left\|W_p(\mu_u,\hat\nu_u)\right \|_p \leq C_2\text{,}
\end{equation*}
where $p\geq 2$ is the same as that in Hypothesis (\ref{hyp1}).

\end{theorem}

{\bf Remark:} Note that under Hypothesis (\ref{hyp1}), the propagation of chaos holds, which will be used in the proof.  The three items (\ref{cod1})-(\ref{cod3}) can be computed and controlled by optimizing the corresponding neural network, where $\epsilon$ can be controlled by setting a threshold for the neural network loss. When $N \to \infty$ and $\epsilon_1, \epsilon_2, \epsilon_3 \to 0 $, $C_1$ and $C_2$ will converge to $0$.
\begin{proof}
    We leave the details of this proof in the Appendix
\end{proof}

Next we present the error analysis of the approach by using the same neural network $F$. Let $Y_t:= F_{\theta}(t,W_t)$, $Y_t^n:= F_{\theta}(t,W_t^n)$ where  $\{W^n, 1\le n\le N\}$ are $N$ independent Brownian motions, and $\hat{\nu}_t^{(N)}:=\frac{1}{N} \sum_{n=1}^{N}\delta_{Y^{n,N}_t}$. Note that in this case $Y_t^n, 1\le n\le N$ are independent and the main difference from the previous result lies in the fact that there is no need to use propagation of chaos.

\begin{theorem}\label{theorem2}
 Consider the MV-SDEs \ref{eq1} and assume the Hypothesis (\ref{hyp1}). Denote the training time points by $t_0=0,t_1, \dots ,t_m=T \in [0,T]$. Assume that there exist three constants $\epsilon_1,\epsilon_2,\epsilon_3>0$ such that the following conditions hold:
\begin{eqnarray}\label{codd1}
\|X_0 - Y_0 \|_p \leq \epsilon_1   \text{,}\\
\label{codd2}\int_{0}^{t} \left\|b(s,Y_s, \hat\nu_s^{(N)}) -  b^{Pseudo}_{n,s} \right \|_p ds   \leq \epsilon_2 \text{,}\\
\label{codd3}\left\| \sqrt{\int_{0}^{t} \left\| \sigma(s,Y_s, \hat\nu_s^{(N)}) -  \sigma^{Pseudo}_{n,s} \right\|^2 ds}\right\|_p \leq \epsilon_3 \text{.}
\end{eqnarray}
Then there exist two constants $C_1$ and $C_2$ depending on $L,T,N,d,\epsilon_1,\epsilon_2,\epsilon_3$, such that
\begin{align*}
\left \| \sup _{u \in [0,t]}|X_u-Y_u|
\right \|_p\leq
 C_1 \nonumber
\end{align*}
and
\begin{align*}
\sup_{u \in [0,t]}& \left\| W_p(\mu_u, \hat{\nu}_u) \right\|_p \leq  C_2,
\end{align*}
where $p\geq 2$ is the same as that in Hypothesis (\ref{hyp1}).
Furthermore,  $C_1$ and $C_2$ will converge to $0$ when $N \to \infty$ and $\epsilon_1, \epsilon_2, \epsilon_3 \to 0 $.
\end{theorem}
\begin{proof}
    We leave the details of this proof in the Appendix. Note that the rate of convergence is different in this case, and this convergence is slower compared to the Theorem 3.6.
\end{proof}

Although this proof is similar to the previous one, it is worth noting that the approximating goal here is the MV-SDE instead of the particle system (\ref{eq2}).  Now there are many methods using supervised learning to approximate the MV-SDE, such as \cite{han2024learning}. The distribution $\nu_t$ of $Y_t$ is determined by $F$, but it is almost impossible to obtain the closed form. Thus for practical computations, we approximate $\nu_t$ using the empirical distribution $\hat{\nu_t}$ of i.i.d. samples, based on Monte Carlo integration rather than relying on the propagation of chaos for interacting particle system. So we do not assume the propagation of chaos in advance for this case, and we use the result in \cite{Fournier2013OnTR} for the \textit{i.i.d.} case.

Next, we perform the analysis for the self-interaction case. We first state the following hypothesis under which the MV-SDE has the unique stationary distribution, see \cite{du2023empirical},.
\begin{hypothesis}\label{hyp2}
The coefficients $b$  and  $\sigma$ are continuous functions on  $\mathbb{R}^{d} \times \mathcal{P}_2\left(\mathbb{R}^{d}\right)$  and bounded on any bounded set; there are constants  $\alpha>\beta \geq 0$, $\gamma>0$, $\rho>0 $, and  $K \geq 0 $ such that
\begin{eqnarray*}
\begin{aligned}
2\langle b(x, \mu)-b(y, \nu), x-y\rangle+\|\sigma(x, \mu)-\sigma(y, \nu)\|^{2} \leq-\alpha|x-y|^{2}+\beta W_{2}(\mu, \nu)^{2} \\
2\langle b(x, \mu), x\rangle+(1+\rho)\|\sigma(x, \mu)\|^{2} \leq-\gamma|x|^{2}+K\left[1+\rho+\mu\left(|\cdot|^{2}\right)\right] \text{,}
\end{aligned}
\end{eqnarray*}
for all  $x, y \in \mathbb{R}^{d}$  and  $\mu, \nu \in \mathcal{P}_2\left(\mathbb{R}^{d}\right)$.
\end{hypothesis}

\cref{hyp2} is a common Hypothesis to make MV-SDE have a unique stationary distribution. When we consider the algorithms to simulate the equation \cref{eq9}, we also need a \cref{hyp1}( $b$ and $\sigma$ without time $t$) to analyze the error of MV-SDE. Under the \cref{hyp1} and the \cref{hyp2}, we have the following estimation of the moment and Wasserstein distance bound in the self-interaction particle system.

\begin{theorem}[Estimation of moments and Wasserstein distance bound in self-interaction particle system] Assume that the Hypothesis \cref{hyp1} and \cref{hyp2} hold. Let $X_t$ be the solution of \cref{eq9} and $Y_t = F(t,W_t)$ with $F$ a neural network. 
Let $ b ^{Pseudo}_{t} := F_1(t,W_t)+ \frac{1}{2}F_{22}(t,W_t) $ and
$ \sigma ^{Pseudo}_{t} :=F_2(t,W_t) $. Given $N\ge 1$, let $t_k=\frac{k}{N}T, k=0,1\cdot, N$. For any $t\in (0,T]$ with $t_k\le t< t_{k+1}$, let $\hat\kappa_{t}= \frac{1}{t_k} \sum_{i=1}^{k} \delta_{X_{t_i}}$ and $\hat\varrho_{t} = \frac{1}{t_k} \sum_{i=1}^{k} \delta_{Y_{t_i}}$. Assume that there exist three constants $\epsilon_1,\epsilon_2,\epsilon_3>0$ such that the following conditions hold:
\begin{eqnarray} \label{cod4}
\|X_0 - Y_0 \|_2 \leq \epsilon_1  \text{,}\\
\label{cod5}
\int_{0}^{t} \left \|b(Y_t, \hat \varrho_{t}) - b ^{Pseudo}_{s}\right \|_2 ds   \leq \epsilon_2 \text{,}\\
\label{cod6}
\left  \| \sqrt{\int_{0}^{t} \left  \| \sigma(Y_s, \hat\varrho_{s}) - \sigma ^{Pseudo}_{s}\right\|^2 ds}\right\|_2 \leq \epsilon_3 \text{.}
\end{eqnarray}
Then there exist two constants $C_1$ and $C_2$ depending on $L,T,d,\epsilon_1,\epsilon_2,\epsilon_3$, such that
\begin{eqnarray}
\left \| \sup _{u \in [0,t]}|X_u-Y_u|
\right \|_2\leq
 C_1\text{,} \\
 \sup _{u  \in [0,t]}\left\|W_2(\hat\kappa_{u},\hat\varrho_{u})\right \|_2
\leq
 C_2\text{.}
\end{eqnarray}
Furthermore, $C_1$ and $C_2$ will converge to $0$ when $T \to \infty$ and $\epsilon_1, \epsilon_2, \epsilon_3 \to 0 $.
 \begin{proof}
     In the Appendix
 \end{proof}

\end{theorem}

\textbf{Note: }
In the discussion of neural network approximation, in our problem, we do not consider the generalization issue. Specifically, we focus only on the solution under a given known Brownian motion and do not aim to use the obtained architecture $F$ for prediction. Therefore, we referenced a classical result \cite{Ryck2021OnTA}, which states that neural networks with the tanh activation function can approximate functions with second-order derivatives. The activation function with tanh to construct $Y^{n,N}_t$, since it can be shown that for any function $g \in W^{2,p} $, there exists a tanh neural network $\hat g \in \mathcal{NN}^n_m $ such that $||g-\hat g||_{L^{\infty}}<\mathcal{\varepsilon}_{n,m,d}$. This detail has been omitted in our discussion.

\section{Algorithms and simulations}
\subsection{Algorithms}
For McKean-Vlasov equation with determined initial values we can often control the following terms (\ref{cod1}, \ref{cod2}, \ref{cod3} or \ref{cod4}, \ref{cod5}, \ref{cod6}) to obtain the corresponding algorithms.

\begin{algorithm}
\caption{DPM for McKean-Vlasov equation
in particle system}
\begin{algorithmic}
\STATE{Define: epochs: $K$, Total point in time: $D$, Learning Rate: $r$, Initial value :$X_0$, Brownian motion :$W_t^n$, Time Series: $t_0,t_2,\dots,t_D=T $.\newline
Neural network: $Y_t^{n,N} = F^n(t,w) = F_{\theta_n}^n(t,w)\text{,} \quad n \in \left\{1,2,\dots,N\right\} $ and $\theta_n$ is the parameter of a neural network. $\lambda_i$ is the weights of the loss function. $\varepsilon = min(\varepsilon_1,\varepsilon_2,\varepsilon_3)$ is the required error threshold.}
\STATE{Note: MSE($a$,$b$) is the value to calculate the mean square error loss between $a$ and $b$. If using a non-interactive particle system here $F^n(t,w)$ can be modified to $F(t,w)$.

}
\FOR{$k$ in $1:K$}

\STATE {Calculate loss}
\begin{align*}
&\mathcal{L}_1= \sum_{i=1}^{N} MSE(X_0,F^i(0,W_0^n)). \\
&\mathcal{L}_2= \sum_{i=1}^{N} \sum_{t=1}^{D}\lambda_iMSE(
b(t,F_{\theta_i}^{i,N}(t,W_t^n),\hat\nu_t),
\frac{\partial F_{\theta_i}^{i,N}(t,W_t^n)}{\partial t}+\frac{1}{2}
\frac{\partial^2 F_{\theta_i}^{i,N}(t,W_t^n)}{\partial w^2}  ). \\
&\mathcal{L}_3= \sum_{i=1}^{N} \sum_{t=1}^{D}\lambda_i MSE(
\sigma(t,F_{\theta_i}^{i,N}(t,W_t^n),\hat\nu_t),
\frac{\partial F_{\theta_i}^{i,N}(t,W_t^n)}{\partial w}).
\end{align*}

\FOR{$n$ in $1:N$}
\STATE{Update parameters $\theta_n^i \leftarrow \theta_n^{i-1} - \nabla_\theta (\mathcal{L}_1+\mathcal{L}_2+\mathcal{L}_3)*r $}.
\ENDFOR
\IF{$ \mathcal{L}_1<\varepsilon_1$ and
$ \mathcal{L}_{2}<\varepsilon_2$ and $ \mathcal{L}_3<\varepsilon_3$ }
\STATE{Break}
\ENDIF
\ENDFOR
\RETURN {$Y_t^{n,N}$}
\end{algorithmic}
\end{algorithm}

\begin{algorithm}
\caption{DPM for McKean-Vlasov equation in self-interaction particle system}
\begin{algorithmic}
\STATE{Define: epochs: $K$, Total point in time: $D$, Learning Rate: $r$, Initial value :$X_0$, Brownian motion :$W_t$, Time Series: $t_0,t_2,\dots,t_D=T $.\newline
Neural network: $Y_t = F(t,w) = F_{\theta}(t,w) $ and $\theta$ is the parameter of a neural network. $\lambda_i$ is the weights of the loss function. $\varepsilon = min(\varepsilon_1,\varepsilon_2,\varepsilon_3)$ is the required error threshold.}
\STATE{Note: MSE($a$,$b$) is the value to calculate the mean square error loss between $a$ and $b$.}

\FOR{$k$ in $1:K$}

\STATE {Calculate loss}
\begin{align*}
&\mathcal{L}_1=  MSE(X_0,F(0,W_0))\text{,}\\
&\mathcal{L}_2= \sum_{t=1}^{D}\lambda_iMSE(
b(F_{\theta}(t,W_t^n),\hat\kappa_{t}),
\frac{\partial F_{\theta}(t,W_t)}{\partial t}+\frac{1}{2}
\frac{\partial^2 F_{\theta}(t,W_t)}{\partial w^2}  )\text{,} \\
&\mathcal{L}_3= \sum_{t=1}^{D}\lambda_i MSE(
\sigma(F_{\theta}(t,W_t),\hat\kappa_{t}),
\frac{\partial F_{\theta}(t,W_t)}{\partial w})\text{.}
\end{align*}
\STATE{Update parameters $\theta^i \leftarrow \theta^{i-1} - \nabla_\theta (\mathcal{L}_1+\mathcal{L}_2+\mathcal{L}_3)*r $}
\IF{$ \mathcal{L}_1<\varepsilon_1$ and
$ \mathcal{L}_{2}<\varepsilon_2$ and $ \mathcal{L}_3<\varepsilon_3$ }
\STATE{Break}
\ENDIF
\ENDFOR
\RETURN {$Y_t$}
\end{algorithmic}
\end{algorithm}

\subsection{Simulations}
\subsubsection{Compare the order of convergence of the two methods}
\begin{figure}[h]
    \centering
    \includegraphics[height=1.9in, width=3.0in]{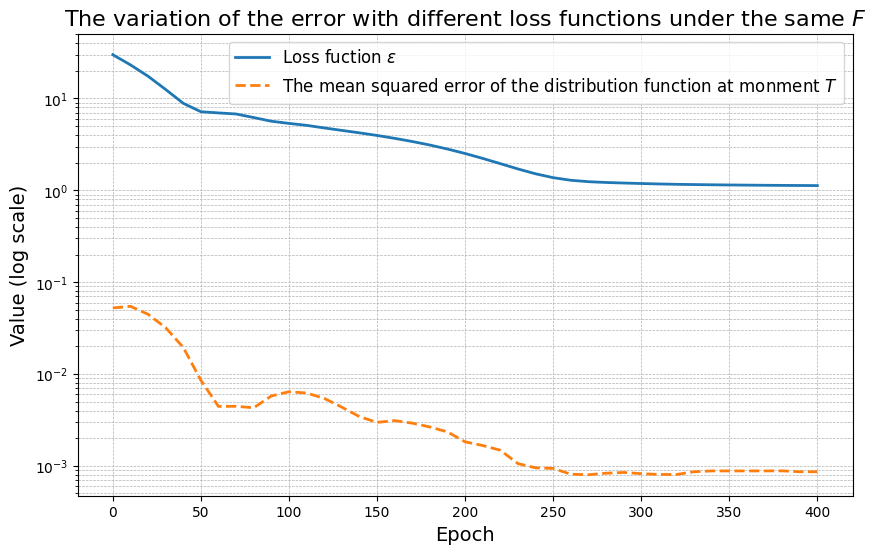}
    \caption{loss functions vs $MSE_{dist}$  (Same $F$)}
    \label{fig:enter-label1}
\end{figure}

\begin{figure}[h]
    \centering
    \includegraphics[height=1.9in, width=3.0in]{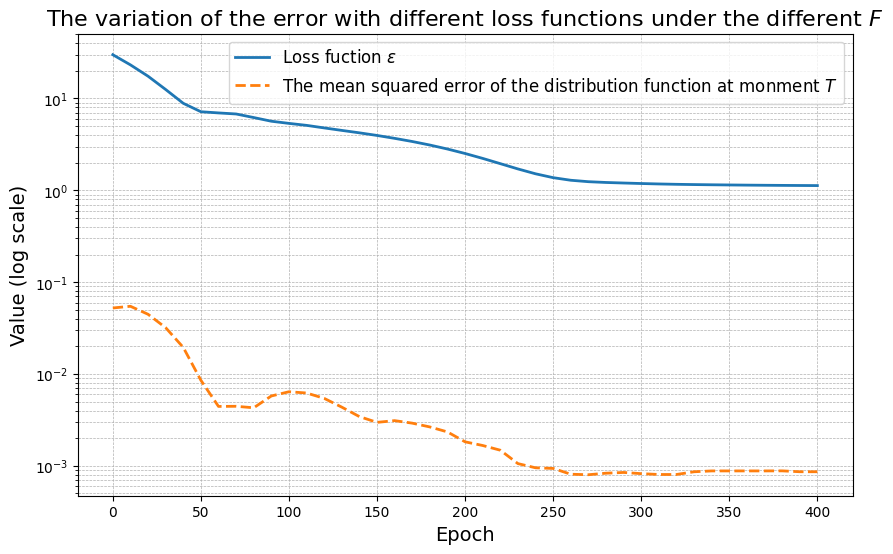}
    \caption{loss functions vs $MSE_{dist}$ (Different $F$)}
    \label{fig:enter-label2}
\end{figure}

\begin{figure}[h]
    \centering
    \includegraphics[height=1.9in, width=3.0in]{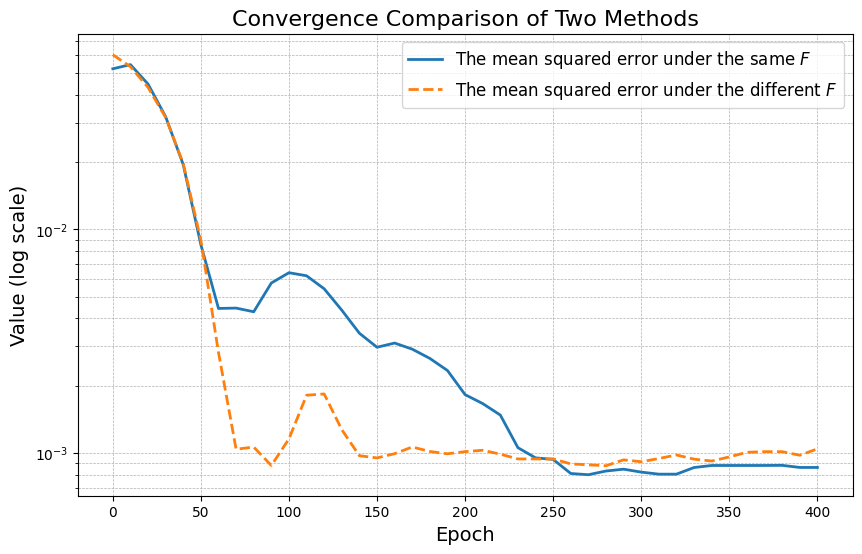}
    \caption{Comparison of $MSE_{dist}$}
    \label{fig:convergence}
\end{figure}

In our previous statements, we proposal two versions of algorithms for the particle system, which aim to approximate \textbf{Auxiliary system} (\ref{eq3}) and \textbf{Interacting particle system} (\ref{eq2}), leading to two different theoretical results regarding the order of convergence. We begin our experiments with a simple example, i.e., the Burgers' equation, for more precise in \cite{Bossy1997ASP,Liu2022ParticleMA}.
\[
d X_{t}=\int_{\mathbb{R}}^{} H(X_t-y)d\mu_t(dy)dt+\sigma dW_t \quad H(x) = \mathbb{I}_{(x\geq 0)}\text{.}
\]
The cumulative distribution function of the $X_1$ has the explicit structure when $X_0 = 0$:
\[
F(x) = \frac{\int_{\mathbb{R}}^{+} \exp \left( -\frac{1}{\sigma^2}[\frac{(x-y)^2}{2}+y] \right) dy}{\int_{\mathbb{R}}^{+} \exp \left( -\frac{1}{\sigma^2}[\frac{(x-y)^2}{2}+y\mathbb{I}_{(y\geq0)}] \right) dy} \text{,}
\] where \(\mathbb{I}_{(y \geq 0)}\) is the indicator function. We choose $h = 0.0001$ and $N= 100000$ to obtain the "true" solution of Auxiliary system by using equation (\ref{eq4}) denotes as PEM (Particle Euler-Maruyama). First, we run the algorithms separately and record the loss function and the mean square error of the distribution at each iteration, which is given by: 
\begin{equation*}
  MSE_{dist} :=\frac{1}{k}\sum _k |F_{Y_T}(x_k) - F_{X_T}(x_k) |^2
\end{equation*}
where $x_k$ is sampled from a uniform distribution in $(-1.5,2.5)$. $F_X$ refers to the distribution function of $X$. We set \( h = 0.01 \), \( M = 101 \), and \( N = 50 \), which represent the time interval of Brownian motion, the number of time points, and the number of trajectories, respectively. We compare the relationship between the loss function and the error in two different cases in Figures \ref{fig:enter-label1},\ref{fig:enter-label2}. Furthermore, the convergence results of the two methods are shown in Figure \ref{fig:convergence}.

We compare these two methods at the path level. For comparison, we simulated the trajectory under this parameter setting using PEM. In our theoretical framework, the approximate results obtained by using different \( F \) should align with those obtained using the interacting particle system (\ref{eq2}). However, the results obtained with the same \( F \) will be somewhat less accurate, as the particles are independent in this case.

\begin{figure}
    \centering
    \includegraphics[height=1.9in, width=3.0in]{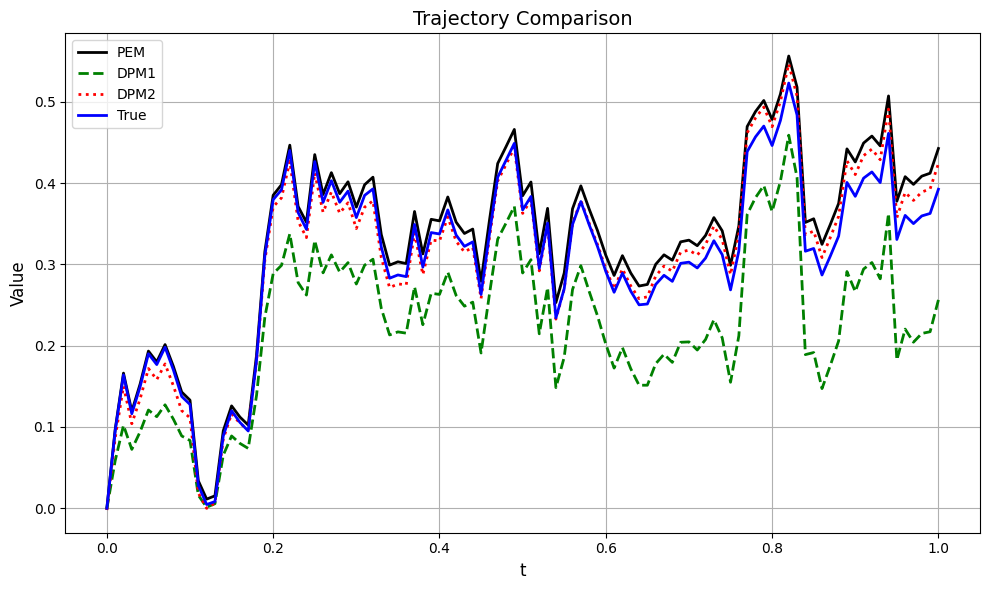}
    \caption{Comparison of the solution trajectories under different methods.}
    \label{fig:trajectories}
\end{figure}

This can be observed from Figure \ref{fig:trajectories}. DPM1, which uses the same \( F \) without assuming chaotic propagation, converges more slowly. On the other hand, DPM2 converges to the classical PEM solution, meaning the \( \text{loss} \) reaches 0, which is equivalent to \( h \to 0 \) in PEM. However, the error between these methods and the true solution still contains a component caused by \( N \).

For a more detailed comparison, we calculate the error at each time step and average it over the number of trajectories. The trajectory error is defined as follows:
\begin{equation*}
\text{Trajectory error} = \frac{1}{N} \sum_{i=1}^{N} |Y_t^i - X_t^{\text{true},i}|^2 
\end{equation*}

\begin{figure}
    \centering
    \includegraphics[height=1.9in, width=3.0in]{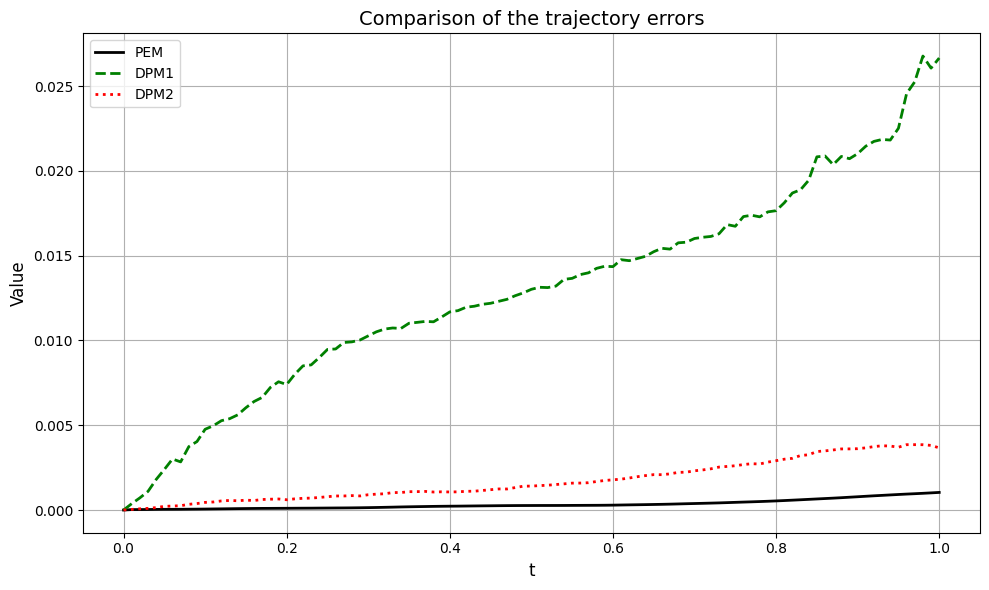}
    \caption{Comparison of the trajectory errors.}
    \label{fig:Trajectoryerror}
\end{figure}
The final results are shown in Figure \ref{fig:Trajectoryerror}. From the figure, we can see that these methods increase as \( t \) grows, which is due to the constant in front of the error term growing exponentially with time. That DPM2 and the PEM do not completely coincide is due to errors introduced by neural network optimization and approximation. These errors can be reduced through proper tuning of hyperparameters and adjustments to the neural network architecture. 
.

\textbf{Compare with the PEM method}.
We compute the cumulative distribution functions of the conventional particle method and the deep learning-based particle method and compare them with the Euler-Maruyama method (PEM). In the experiment, we chose the point with equal step size (which can be non-equal step size in fact) for calculation, and we simulate the trajectory of $N$ independent Brownian motion.
We set $\sigma = 0.5$, $h= 1$, $N = 50$, $M = 2$, where $h$ is the step size, $M-1$ is the number of iterative steps, and $N$ is the number of particles used. We use a 2-layer neural network with 32 units per layer with activation function \textit{tanh()}, and each particle is trained for 150 rounds using a separate neural network with a separate optimizer Adam with a learning rate of 0.001, and we scale the distribution function to the interval $[-1.5,2.5]$, obtain the following image and calculate the $MSE_{dist}$ of the two methods for $F(x)$. Here, we focus more on accuracy, which is why we used DPM2, i.e. the neural network under different $F$.
\begin{figure}[h]
\centering
\includegraphics[height=1.9in, width=3.0in]{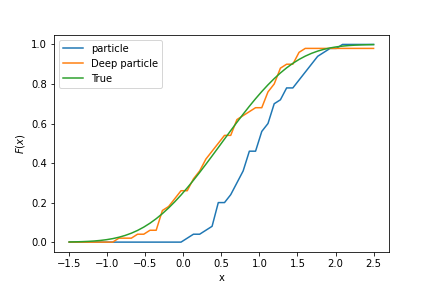}
\caption{Comparison using the particle method and the deep particle method $F(x)$}
\label{fig:particle_vs_deep_particle}
\end{figure}

\begin{figure}[h]
\centering
\includegraphics[height=1.9in, width=3.0in]{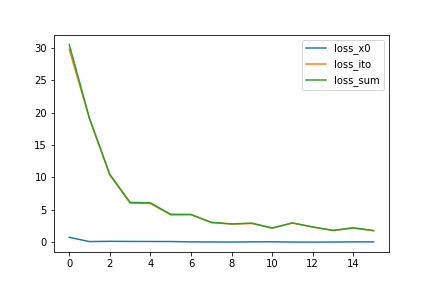}
\caption{Loss function for training}
\label{fig:loss_function_training}
\end{figure}

In this extreme condition, the step size is 1 and there are only 2 iterations. We achieve $MSE= 0.0005959$, while the particle method based on the Euler iteration has $MSE= 0.02812570$. We compare the $MSE$ obtained with two different methods for different numbers of paths at different step sizes as Table \ref{1}.

\begin{table}
\label{1}
\centering
\begin{tabular}{|l|l|l|l|l|l|l|}
\hline
Method & $h$ & $M$  & $N$&  $Mean Square Error$ & Training time(s)\\ \hline

PEM& 1 & 2 &  50 & 2.812570e-02 & $\times$ \\ \hline

DPM  & 1 & 2 &  50 & 5.959000e-04 & 53.51  \\ \hline

PEM& 0.5 & 3 &  50 & 6.821674e-03 & $\times$ \\ \hline

DPM  & 0.5 & 3 &  50 &  4.813673e-04 & 54.93  \\ \hline

PEM& 1 & 2 &  100 & 2.486127e-02 & $\times$ \\ \hline

DPM  & 1 & 2 &  100 & 1.155264e-03 & 107.73  \\ \hline

PEM& 0.5 & 3 &  100 & 5.120935e-03 & $\times$ \\ \hline

DPM  & 0.5 & 3 &  100 &5.906984e-04 & 110.31\\ \hline
\end{tabular}
\caption {Comparison using the particle method and the deep particle method }
\end{table}
The above experiments are completed under the unified experimental configuration, and the epoch is 100, and learning rate is 0.001. In order to better compare the performance of the models with the same parameter configuration, we did not adjust the experiments for the larger $N$. An important observation is that, for the PEM, the algorithm's running time increases rapidly as the number of iterations grows and number of trajectories. In contrast, in DPM, this increase only results in more points being added to the training set. Since the training of the neural network can be parallelized, this does not lead to an unacceptable increase in time. However, increasing the number of training points will increase the model's training time for a given level of accuracy. With sufficient GPU memory, this makes it possible to simulate large-scale systems.

\subsubsection{Simulate particle systems in the case of meshless}
The next experiment demonstrates a feature that PEM does not possess. We can solve the MV-SDE in the absence of the Brownian motion trajectories. The missing Brownian motion may cause a large step size at some point, which could lead to the solution exploding. In the following experiment, we demonstrate the model's performance with unequal step sizes, a feature not supported by traditional solvers. To represent the unequal step sizes, we randomly delete the paths of the sampled Brownian motion using $numpy.random.choice()$ and we record the performance of the model at different deletion ratios.

The censoring rule is to randomly censor but keep the initial and final points, which is the closest to the real situation, since the noisy data we observe may be unequally spaced, or have a large proportion of missing values.
\begin{figure}[h]
\centering
\includegraphics[height=1.9in, width=3.0in]{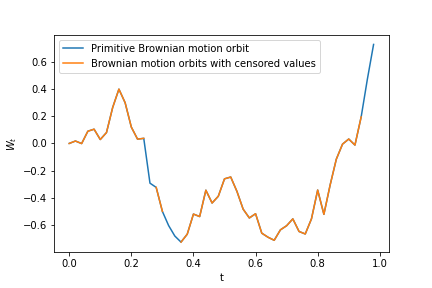}
\caption{Deletion rate is 0.2}
\label{fig:deletion_rate_0.2}
\end{figure}

\begin{figure}[h]
\centering
\includegraphics[height=1.9in, width=3.0in]{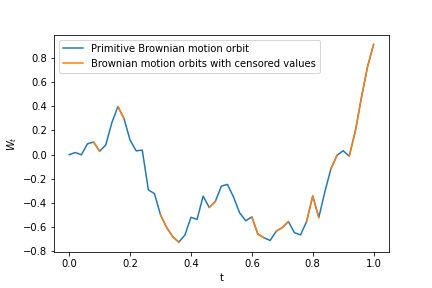}
\caption{Deletion rate is 0.5}
\label{fig:deletion_rate_0.5}
\end{figure}

\begin{table}
\centering
\begin{tabular}{|l|l|l|l|l|l|l|}
\hline
Method & $h$ & $M$  & $N$  & Deletion rate&  Mean Square Error &  time(s)\\ \hline

DPM  & 0.02 & 51 &  50 & 0.1 &  4.948411e-04 & 129.50 \\ \hline

DPM  & 0.02 & 51 &  50 & 0.2 &  5.023381e-04 & 117.58 \\ \hline

DPM  & 0.02 & 51 &  50 & 0.5 &  6.261634e-04 & 94.29 \\ \hline

DPM  & 0.02 & 51 &  50 & 0.7 &  5.065099e-04 & 80.00  \\ \hline

DPM  & 0.02 & 51 &  50 & 0.9 &  7.486928e-04 & 64.22  \\ \hline

DPM  & 0.01 & 101 &  50 & 0.1 &  1.050783e-03 & 201.13 \\ \hline

DPM  & 0.01 & 101 &  50 & 0.2 &  1.103091e-03 & 187.14\\ \hline

DPM  & 0.01 & 101 &  50 & 0.5 &  9.647387e-04 & 130.09\\ \hline

DPM  & 0.01 & 101 &  50 & 0.7 &  1.058312e-03 & 97.09\\ \hline

DPM  & 0.01 & 101 &  50 & 0.9 &  9.555768e-04 & 63.90\\ \hline
\end{tabular}
\caption {Comparison using the depth particle method at different deletion rates}
\end{table}

It can be seen that the step size and its deletion rate do not have a significant effect on the final result. This is because this method is dependent on the optimisation process of the neural network, and the construction of the neural network, while the data points are only used as training data.

\subsubsection{Simulate McKean-Vlasov Stochastic Differential Equation for other noises}
We extend this idea to the McKean-Vlasov Stochastic Differential driven by Fractional Brownian motion. It is defined as follows:
\begin{equation*}
\begin{aligned}
\left\{\begin{array}{l}
d X_{t}=b\left(t, X_{t}, \mu_{t}\right) d t+\sigma\left(t, X_{t}, \mu_{t}\right) d W_{t}^H \\
\mu_{t} = \mathcal{L}(X_{t}) \text{ is distribution of $X_t$}
\\
X_0=\mathcal{L}_0 \text{ is distribution of $X_0$} \text{,}\quad\quad t \in [0, T]
\end{array}\right.\text{,}
\end{aligned}
\end{equation*}
$W_{t}^H$ is the fractional Brownian motion with Hurst is $H \in [\frac{1}{2},1)$.
In this case we need to modify the stochastic chain law in the algorithm:
\begin{eqnarray*}
Y_t^{n,N} &=& F^n(t,W_t^{n,H}) \text{,} \\
dY_t^{n,N} &=& F^n_1(t,W_t^{n,H})dt + Ht^{2H-1}F^n_{22}(t,W_t^{n,H})dt + F^n_2(t,W_t^{n,H})dW_t^{n,H} \text{.}
\end{eqnarray*}
The same compatibility can also be demonstrated using a similar approach \cite{he2022explicit}.
For fractional Brownian motion we use the method in \cite{Muniandy2001ModelingOL}.
\begin{eqnarray*}
dX_t = \int_{\mathbb{R}}(X_t- y )\mu_t(dy) + \sigma dW_t^H\text{,} \quad X_0 = x_0\text{,}
\end{eqnarray*}
Our setup parameter is $\sigma = 0.5$, $ H =0.7$, $x_0 = 0$. Since this equation has no analytical solution, we use a smaller step size ($h = 0.05$) and a larger number of paths $N = 500$ to obtain the ''true'' solution.
Similarly we compute the cumulative distribution function of the solution at the moment $T= 1$.
\begin{figure}[h]
\centering
\includegraphics[height=1.9in, width=3.0in]{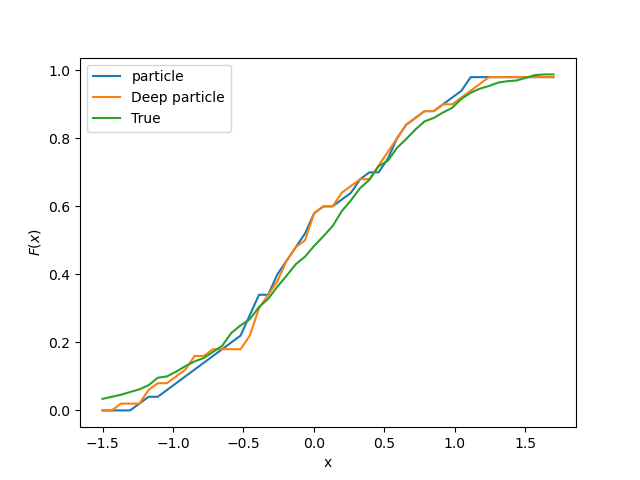}
\caption{Comparison using the particle method and the deep particle method $F(x)$.}
\label{fig:compare}
\end{figure}

\begin{figure}[h]
\centering
\includegraphics[height=1.9in, width=3.0in]{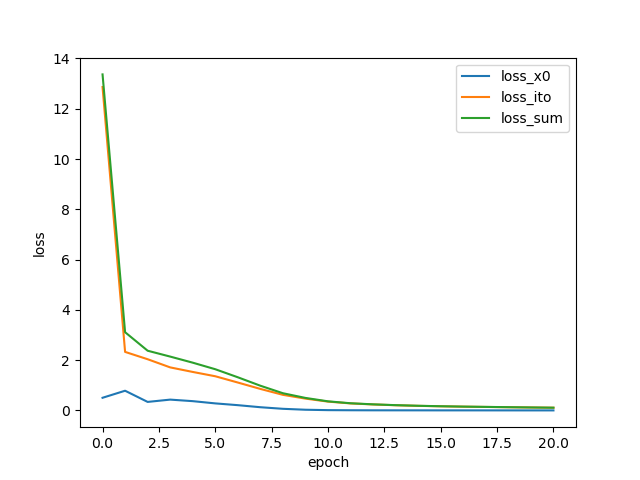}
\caption{Loss function for training.}
\label{fig:loss}
\end{figure}

The particle method is obtained by using the standard Euler simulation, in which we use $h = 0.1$, $M= 11$, in the DPM we set the neural network with \textit{tanh()}, the size of the hidden layer is 64, and 2 hidden layers are used. In the training we recorded the value of the loss function every 10 times and trained 200 epochs with the Adam optimizer with a learning rate of 0.001. With this parameter configuration, the Euler method is used to obtain MSE is 1.345920e-03, the MSE obtained using the DPM method is 1.165120e-03. Further we consider the results of this approach for different parameter configurations. We just choose one type of noise. But this method has great flexibility and can cover most of the MV-SDE.

\begin{table}[h]
\centering
\begin{tabular}{|l|l|l|l|l|l|l|}
\hline
Method & $h$ & $M$  & $N$&  $Mean Square Error$ & Training time(s)\\ \hline

PEM & 0.1 & 11 &  50 & 1.345920e-03 & $\times$ \\ \hline
DPM  & 0.1 & 11 &  50 &  1.165120e-03 & 80.13\\ \hline

PEM & 0.5 & 3 &  50 & 6.825920e-03 & $\times$ \\ \hline
DPM  & 0.5 & 3 &  50 &  4.995520e-03 & 68.72  \\ \hline

PEM & 0.1 & 11 &  100 & 1.582320e-03 & $\times$ \\ \hline
DPM  & 0.1 & 11 &  100 &  2.713520e-03 & 161.84\\ \hline

PEM & 0.5 & 3 &  100 & 2.823920e-03 & $\times$ \\ \hline
DPM  & 0.5 & 3 &  100 & 1.663920e-03  &  142.35\\ \hline

\end{tabular}
\caption {Comparison using the particle method and the deep particle method}
\end{table}

\subsubsection{Simulate stationary distributions for self-interacting progress}

\begin{figure}[h]
\centering
\includegraphics[height=1.8in, width=2.5in]{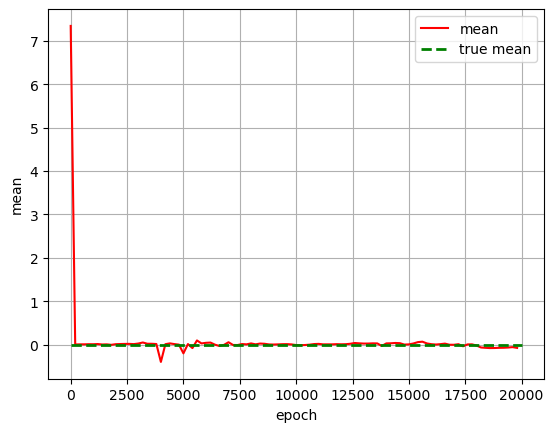}
\caption{The mean value of the solution $X_t$ during training.}
\label{fig:The mean value of the solution $X_t$ during training}
\end{figure}

\begin{figure}[h]
\centering
\includegraphics[height=1.8in, width=2.5in]{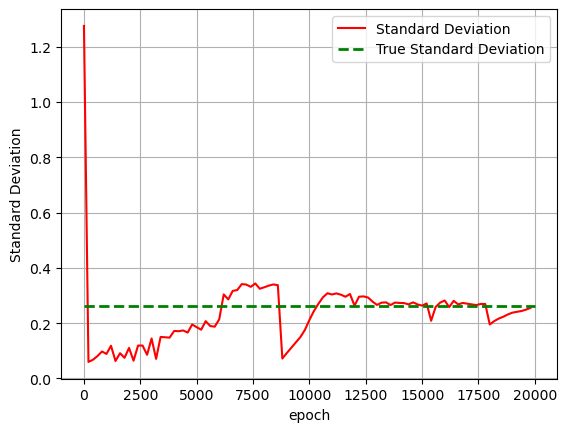}
\caption{The standard deviation value of the solution $X_t$ during training.}
\label{fig:The standard deviation value of the solution $X_t$ during training.}
\end{figure}

\begin{figure}[h]
\centering
\includegraphics[height=1.8in, width=2.5in]{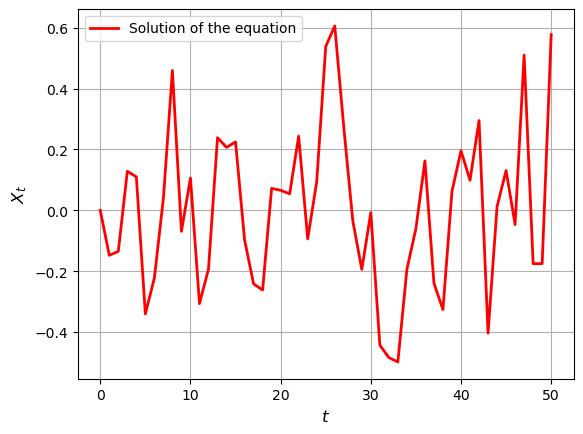}
\caption{The path of the solution after training is complete.}
\label{fig:The path of the solution after training is complete.}
\end{figure}

We consider the following MV-SDE that satisfies \cref{hyp1} and \cref{hyp2}. We use one path of Brownian motion to approximate the stationary distribution of this equation by using the self-interacting DPM method.
\begin{eqnarray*}
dX_t = (-2X_t - \mathbb{E}X_t )dt + \sigma dW_t \text{,}\quad X_0 = x_0\text{,}
\end{eqnarray*}
where $\sigma = 0.5$ and $X_0 = 0$. At first, we use Euler-Maruyama combined with the interacting particle system to obtain the explicit solution. The stationary distribution of this equation is a normal distribution, and it converges to this stationary distribution for $T\geq1$. Euler Maruyama's parameters are set to $N = 1000$ and the step size $h = 0.01$ and the number of steps is $M = 1001$. The approximation we get is that the mean of a stationary distribution is -0.0002 and the standard deviation is 0.2633. We then set $M=51$ and we computed the solution for the next 50 steps and trained 20,000 epochs, we recorded the mean and standard deviation of the solution obtained each time. Figures \ref{fig:The mean value of the solution $X_t$ during training} and \ref{fig:The standard deviation value of the solution $X_t$ during training.} show the errors as training progresses, while Figure \ref{fig:The path of the solution after training is complete.} shows the trajectories of the self-interacting solution after training.

It can be seen that this method converges quickly and does not require storing historical record points. For general self-interacting solvers, the following iterative relation is often present:
\begin{equation}
X_{t_i} = Solver(X_{t_0},X_{t_1}\dots,X_{t_i},h,B_{t_{i-1}})\text{.}
\end{equation}
In general settings for solving stationary distributions, a large number of time points is often required. However, as time increases, more iteration steps are needed, which requires more points to compute the empirical distribution at the final time step. In contrast, in DPM, these points are predefined (by providing the Brownian motion and then obtaining the output), and all the points are incorporated into the loss function for training in one go. Thus, it is a global control that does not require storing a large number of points. We only need to optimize over the given set of points at each training step. Moreover, since this method is grid-free, the points we provide do not require arbitrarily small time intervals, which is another significant advantage.

\section{Conclusions}
We propose an alternative approach to simulate the McKean-Vlasov equation by combining deep learning to approximate the drift diffusion coefficients. This approach can be considered as a PINN method in the field of stochastic differential equations. This approach eliminates the dependence on the step size and provides error estimates by optimizing the approximation using an appropriate algorithm. The simulations performed support the advantages of this approach. An obvious benefit of this approach is that the discrete points used for function estimation can be inhomogeneous, reducing the need for small step sizes. Furthermore, we suggest that this framework can be extended to handle a wider range of noise, including McKean-Vlasov SDEs driven by fractional Brownian motion.

\bibliographystyle{siamplain}
\bibliography{references}

\section{Appendix}

\subsection{Proof for Theorem 3.6}
\begin{proof}
By using Lemmas \cref{lemma1}, \cref{lemma2} and Hypothesis (\ref{hyp1}), we have
\begin{eqnarray}
&&\left\|\sup _{u \in [0,t]}|X_u^{n,N}-Y_u^{n,N}|\right \|_p \leq 
\left\| \sup _{u \in [0,t]}|X_0^{n,N}-Y_0^{n,N}
+\int_{0}^{u}b(s,X_s^{n,N},\hat\mu_s)-b ^{Pseudo}_{n,s}ds \right\|_p\nonumber \\
&+&
\left\|\sup _{u \in [0,t]}|\int_{0}^{u}\sigma(s,X_s^{n,N},\hat\mu_s)-\sigma ^{Pseudo}_{n,s}dB_s |  \right\|_p \nonumber \\
&\leq & \left\|X_0^{n,N}-Y_0^{n,N}\right\|_p +\left\|\sup _{u \in [0,t]}|\int_{0}^{u}b(s,X_s^{n,N},\hat\mu_s)-b ^{Pseudo}_{n,s} +b(s,Y_s^{n,N},\hat\nu_s) -b(s,Y_s^{n,N},\hat\nu_s)ds |\right\|_p\ \nonumber \\
&+&
\left\|\sup _{u \in [0,t]}\int_{0}^{u}\sigma(s,X_s^{n,N},\hat\mu_s)-\sigma ^{Pseudo}_{n,s}+\sigma(s,Y_s^{n,N},\hat\nu_s)-\sigma(s,Y_s^{n,N},\hat\nu_s)dB_s |  \right\|_p \nonumber \\ 
&\leq & \left\|X_0^{n,N}-Y_0^{n,N}\right\|_p + \left\|\sup _{u \in [0,t]}\int_{0}^{u}b(s,X_s^{n,N},\hat\mu_t)- b(s,Y_s^{n,N},\hat\nu_s)ds \right\|_p\ \nonumber \\
&+&\left\|\sup _{u \in [0,t]}\int_{0}^{u}b ^{Pseudo}_{n,s}- b(s,Y_s^{n,N},\hat\nu_s)ds \right\|_p
+\left\|\sup _{u \in [0,t]}\int_{0}^{u}\sigma(s,X_s^{n,N},\hat\mu_s)- \sigma(s,Y_s^{n,N},\hat\nu_s)dB_s \right\|_p \nonumber \\
&+&\left\|\sup _{u \in [0,t]}\int_{0}^{u}\sigma ^{Pseudo}_{n,s}- \sigma(s,Y_s^{n,N},\hat\nu_s)dB_s \right\|_p \nonumber \\
&\leq & \epsilon_1 + \epsilon_2 + \epsilon_3C_p+
\left\|\sup _{u \in [0,t]}\int_{0}^{u} b(s,X_s^{n,N},\hat\mu_s)- b(s,Y_s^{n,N},\hat\nu_s)ds \right\|_p\  \nonumber \\
& +&
\left\|\sup _{u \in [0,t]}\int_{0}^{u}\sigma(s,X_s^{n,N},\hat\mu_s)- \sigma(s,Y_s^{n,N},\hat\nu_s)dB_s \right\|_p\  \nonumber \\
&\leq & \epsilon_1 + \epsilon_2 + \epsilon_3C_p+ \int_{0}^{t} \left \|b(s,X_s^{n,N},\hat \mu_s ) - b(s,Y_s^{n,N},\hat \nu_s )\right \|_pds
+
C_{d,p} \left  \| \sqrt{\int_{0}^{t} \left  \| \sigma(s,X_s^{n,N},\hat \mu_s ) - \sigma(s,Y_s^{n,N},\hat \nu_s )\right\|^2 ds}\right\|_p  \nonumber 
\end{eqnarray}
\begin{eqnarray}
&\leq & \epsilon_1 + \epsilon_2 + \epsilon_3C_p+
L\int_{0}^{t}(\left\| X_s^{n,N}- Y_s^{n,N}\right\|_p + \left\| W_p(\hat\mu_s,\hat\nu_s)\right\|_p)ds +  C_{d,p}\left  \| \int_{0}^{t} \left  \| \sigma(s,X_s^{n,N},\hat \mu_s) - \sigma(s,Y_s^{n,N},\hat \nu_s )\right\|^2 ds\right\|^\frac{1}{2}_\frac{p}{2} \nonumber \\
&\leq &  \epsilon_1 + \epsilon_2 + \epsilon_3C_p+
L\int_{0}^{t}\left(\| X_s^{n,N}- Y_s^{n,N}\right\|_p + \left\| W_p(\hat\mu_s,\hat\nu_s)\right\|_p)ds +  C_{d,p}\sqrt{\int_{0}^{t}   \left \| \left \|\sigma(s,X_s^{n,N},\hat \mu_s ) - \sigma(s,Y_s^{n,N},\hat \nu_s ) \right\|  \right\|^2_p ds} \nonumber \\
&\leq &  \epsilon_1 + \epsilon_2 + \epsilon_3C_p+
L\int_{0}^{t}(\left\| X_s^{n,N}- Y_s^{n,N}\right\|_p + \left\| W_p(\hat\mu_s,\hat\nu_s)\right\|_p)ds + C_{p,d,L}\left \{\int_{0}^{t}\left\| X_s^{n,N}- Y_s^{n,N}\right\|^2_p + \left\| W_p(\hat\mu_s,\hat\nu_s)\right\|^2_pds\right \}^\frac{1}{2}  \nonumber.
\end{eqnarray}


\begin{eqnarray*}
\text{Let }
 \varphi(t) := \epsilon_1 + \epsilon_2 + \epsilon_3C_p+
L\int_{0}^{t}\left\| W_p(\hat\mu_s,\hat\nu_s)\right\|_pds
+ C_{p,d,L}\left \{\int_{0}^{t}\left\| W_p(\hat\mu_s,\hat\nu_s)\right\|^2_pds\right\}^\frac{1}{2},
\end{eqnarray*}

\begin{eqnarray}
\text{then }
\left\|\sup _{u \in [0,t]}|X_u^{n,N}-Y_u^{n,N}|\right \|_p  \leq
L\int_{0}^{t}\left\| X_s^{n,N}- Y_s^{n,N}\right\|_pds
+C_{p,d,L}\left \{\int_{0}^{t}\left\|X_s^{n,N}- Y_s^{n,N} \right\|^2_pds\right\}^\frac{1}{2} + \varphi(t)\nonumber.
\end{eqnarray}

By \cref{lemma3}, we have
\begin{equation}\label{wq1}
\left \| \sup _{u \in [0,t]}|X_u^{n,N}-Y_u^{n,N}|
\right \|_p\leq 2e^{(2L+C_{p,d,L}^2)t}\varphi(t).
\end{equation}

Next we estimate $\left\|W_p(\mu_s,\hat\nu_s)\right \|_p$. Note that
\begin{equation}\label{wq2}
\sup _{s  \in [0,t]}\left\|W_p(\mu_s,\hat\nu_s)\right \|_p  \leq \sup _{s  \in [0,t]}\left\|W_p(\mu_s,\hat\mu_s)\right \|_p +\sup _{s  \in [0,t]}\left\|W_p(\hat\mu_s,\hat\nu_s)\right \|_p.
\end{equation}
The propagation of chaos is here used for estimating the first term in the right hand side (r.h.s. in short) of (\ref{wq2}), and we do not provide further details here. However, it is worth noting that this term will tend to $0$ as $N\to\infty$.

For the second term in the r.h.s. of (\ref{wq2}), we have
\begin{eqnarray}
\sup _{s  \in [0,t]}\mathbb{E}W_p^p(\hat\mu_s,\hat\nu_s)
&\leq& \mathbb{E} \left [   \frac{1}{N} \sum_{n=1}^{N}\sup_{s\in [0,t]}|X^{n,N}_s-Y^{n,N}_s|^p\right ] = \frac{1}{N} \sum_{n=1}^{N}\|\sup_{s\in [0,t]}|X^{n,N}_s-Y^{n,N}_s|\|^p_p \nonumber \\
&\leq& (2e^{(2L+C_{p,d,L}^2)t}\varphi(t))^p .\nonumber
\end{eqnarray}

By the definition of $\varphi(t)$, we get
\begin{eqnarray}
\sup _{s  \in [0,t]}\left\|W_p(\hat\mu_s,\hat\nu_s)\right \|_p &\leq&
2e^{(2L+C_{p,d,L}^2)t} \left ( \epsilon_1 + \epsilon_2 + \epsilon_3C_p+
L\int_{0}^{t}\left\| W_p(\hat\mu_s,\hat\nu_s)\right\|_pds
+ C_{p,d,L}\left \{\int_{0}^{t}\left\| W_p(\hat\mu_s,\hat\nu_s)\right\|^2_pds\right\}^\frac{1}{2} \right )  \nonumber \\
&\leq& A+B\int_{0}^{t}\sup _{m \in [0,s]}\left\| W_p(\hat\mu_m,\hat\nu_m)\right\|_pds
+ C\left \{\int_{0}^{t}\sup _{m \in [0,s]}\left\| W_p(\hat\mu_m,\hat\nu_m)\right\|^2_pds\right\}^\frac{1}{2}, \nonumber
\end{eqnarray}
\begin{eqnarray*}
\text{where } A = 2e^{(2L+C_{p,d,L}^2)t}(\epsilon_1 + \epsilon_2 + \epsilon_3C_p) ,B = 2e^{(2L+C_{p,d,L}^2)t}L \text{ and } C=2e^{(2L+C_{p,d,L}^2)t}C_{p,d,L}\text{.}
\end{eqnarray*}
By using \cref{lemma3} once again,  we can obtain
\begin{eqnarray}
\sup _{s  \in [0,t]}\left\|W_p(\hat\mu_s,\hat\nu_s)\right \|_p
\leq 2e^{(2B+C^2)t}A
\end{eqnarray}

Taking this estimate back to $\varphi(t)$, and finally we get the desired results according to  \cref{wq1} and \cref{wq2} respectively.

\end{proof}

\subsection{Proof for Theorem 3.7}
\begin{proof}
By using Lemmas \cref{lemma1}, \cref{lemma2} and Hypothesis (\ref{hyp1}), we have
\begin{eqnarray}
&&\left\|\sup _{u \in [0,t]}|X_u-Y_u| \right \|_p  \leq
\left\|X_0-Y_0\right\|_p +\left\| \sup _{u \in [0,t]}|\int_{0}^{u}b(s,X_s,\mu_s)-b^{Pseudo}_{n,s}ds| \right\|_p\nonumber+\left\|\sup _{u \in [0,t]}|\int_{0}^{u}\sigma(s,X_s,\mu_s)-\sigma ^{Pseudo}_{n,s}dB_s |  \right\|_p \nonumber \\
&\leq & \left\|X_0-Y_0\right\|_p +\left\|\sup _{u \in [0,t]}|\int_{0}^{u}b(s,X_s,\mu_s)  -b(s,Y_s,\nu_s) +b(s,Y_s,\nu_s)-b(s,Y_s,\hat\nu_s^{(N)})+
b(s,Y_s, \hat \nu_s^{(N)})-b ^{Pseudo}_{n,s} ds |\right\|_p\  \nonumber \\
&+&
\left\|\sup _{u \in [0,t]}|\int_{0}^{u}\sigma(s,X_s,\mu_s)-\sigma(s,Y_s,\nu_s) +\sigma(s,Y_s,\nu_s)-\sigma(s,Y_s,\hat\nu_s^{(N)})+\sigma(s,Y_s,\hat\nu_s^{(N)})-\sigma ^{Pseudo}_{n,s} dB_s |  \right\|_p \nonumber \\
&\leq & \left\|X_0-Y_0\right\|_p + \left\|\sup _{u \in [0,t]}\int_{0}^{u}|b(s,X_s,\mu_s)- b(s,Y_s,\nu_s)|ds \right\|_p\ + \left\|\sup _{u \in [0,t]}\int_{0}^{u}|b(s,Y_s,\nu_s)- b(s,Y_s,\hat\nu_s^{(N)})|ds \right\|_p \nonumber \\
&+&\left\|\sup _{u \in [0,t]}\int_{0}^{u}|b ^{Pseudo}_{n,s}- b(s,Y_s,\nu_s)|ds \right\|_p
+\left\|\sup _{u \in [0,t]}\int_{0}^{u}|\sigma(s,X_s,\mu_s)- \sigma(s,Y_s,\nu_s)|dB_s \right\|_p \nonumber 
\end{eqnarray}
\begin{eqnarray}
&+& \left\|\sup _{u \in [0,t]}\int_{0}^{u}|\sigma(s,Y_s,\nu_s)- \sigma(s,Y_s,\hat\nu_s^{(N)})|ds \right\|_p + \left\|\sup _{u \in [0,t]}\int_{0}^{u}|\sigma ^{Pseudo}_{n,s}- \sigma(s,Y_s,\nu_s)|dB_s \right\|_p \nonumber \\
&\leq & \epsilon_1 + \epsilon_2 + \epsilon_3C_p + I(N)+
\left\|\sup _{u \in [0,t]}\int_{0}^{u} |b(s,X_s,\mu_s)- b(s,Y_s,\nu_s)|ds \right\|_p\ +C_{d,p} \left  \| \sqrt{\int_{0}^{t} \left  \| \sigma(s,X_s, \mu_s ) - \sigma(s,Y_s, \nu_s )\right\|^2 ds}\right\|_p   \nonumber \\
&\leq & \epsilon_1 + \epsilon_2 + \epsilon_3C_p + I(N)+
L\int_{0}^{t}(\left\| X_s- Y_s\right\|_p + \left\| W_p(\mu_s,\nu_s)\right\|_p)ds+C_{d,p}\left  \| \int_{0}^{t} \left  \| \sigma(s,X_s,\mu_s) - \sigma(s,Y_s,\nu_s )\right\|^2 ds\right\|^\frac{1}{2}_\frac{p}{2} \nonumber \\
&\leq &  \epsilon_1 + \epsilon_2 + \epsilon_3C_p + I(N)+
L\int_{0}^{t}(\left\| X_s- Y_s\right\|_p + \left\| W_p(\mu_s,\nu_s)\right\|_p)ds + C_{d,p}\sqrt{\int_{0}^{t}   \left \| \left \|\sigma(s,X_s,\mu_s ) - \sigma(s,Y_s, \nu_s ) \right\|  \right\|^2_p ds}\nonumber \\
&\leq &  \epsilon_1 + \epsilon_2 + \epsilon_3C_p + I(N)+
L\int_{0}^{t}(\left\| X_s- Y_s\right\|_p + \left\| W_p(\mu_s,\nu_s)\right\|_p)ds+ C_{p,d,L}\left \{\int_{0}^{t}\left\| X_s- Y_s\right\|^2_p + \left\| W_p(\mu_s,\nu_s)\right\|^2_pds\right \}^\frac{1}{2} \text{,}\nonumber
\end{eqnarray}
 \begin{eqnarray*}\text{where }
I(N):=\left\|\sup _{u \in [0,t]}\int_{0}^{u}b(s,Y_s,\nu_s)- b(s,Y_s,\hat\nu_s)ds \right\|_p + \left\|\sup _{u \in [0,t]}\int_{0}^{u}\sigma(s,Y_s,\nu_s)- \sigma(s,Y_s,\hat\nu_s)ds \right\|_p.
\end{eqnarray*}
By the lipschitz property for $b$ and $\sigma$, we get
\begin{eqnarray*}
I(N) &\leq& L \left\|\sup _{u \in [0,t]}\int_{0}^{u}W_p(\nu_s,\hat\nu_s)ds \right\|_p +  L \left\|\sup _{u \in [0,t]}\int_{0}^{u}W_p(\nu_s,\hat\nu_s)dB_s \right\|_p \nonumber  \\
 &\leq& L \int_{0}^{t} ||W_p(\nu_s,\hat\nu_s)||_p ds + + L \left\{ \int_{0}^{t}||W_p(\nu_s,\hat\nu_s)||^2_p ds \right\} ^{\frac{1}{2}}\text{.}
\end{eqnarray*}

\begin{eqnarray*}
\text{Define } 
 \varphi(t) := \epsilon_1 + \epsilon_2 + \epsilon_3C_p+ I(N)+
L\int_{0}^{t}\left\| W_p(\mu_s,\nu_s)\right\|_pds
+ C_{p,d,L}\left \{\int_{0}^{t}\left\| W_p(\mu_s,\nu_s)\right\|^2_pds\right\}^\frac{1}{2} \text{,}
\end{eqnarray*}

\begin{eqnarray}
\text{then } 
\left\|\sup _{u \in [0,t]}|X_u-Y_u|\right \|_p  \leq
L\int_{0}^{t}\left\| X_s- Y_s\right\|_pds
+C_{p,d,L}\left \{\int_{0}^{t}\left\|X_s- Y_s \right\|^2_pds\right\}^\frac{1}{2} + \varphi(t)\text{.}\nonumber
\end{eqnarray}
By \cref{lemma3}, we have
\begin{equation}
\label{target}
\left \| \sup _{u \in [0,t]}|X_u-Y_u|
\right \|_p\leq 2e^{(2L+C_{p,d,L}^2)t}\varphi(t)\text{,}
\end{equation}
and furthermore

\begin{eqnarray}
&\sup _{s  \in [0,t]}&\left\|W_p(\mu_s,\nu_s)\right \|_p  \leq  \left\|\sup _{u \in [0,t]}|X_u-Y_u|\right \|_p  \leq 2e^{(2L+C_{p,d,L}^2)t}\varphi(t)  \nonumber \\
&\leq& 2e^{(2L+C_{p,d,L}^2)t}\left ( \epsilon_1 + \epsilon_2 + \epsilon_3C_p+ I(N)+
L\int_{0}^{t}\left\| W_p(\mu_s,\nu_s)\right\|_pds
+ C_{p,d,L}\left \{\int_{0}^{t}\left\| W_p(\mu_s,\nu_s)\right\|^2_pds\right\}^\frac{1}{2}  \right )\nonumber \\
&\leq& A+B\int_{0}^{t}\sup _{m \in [0,s]}\left\| W_p(\mu_m,\nu_m)\right\|_pds
+ C\left \{\int_{0}^{t}\sup _{m \in [0,s]}\left\| W_p(\mu_m,\nu_m)\right\|^2_pds\right\}^\frac{1}{2}  \text{,} \nonumber
\end{eqnarray}
\begin{eqnarray*}
\text{where } A = 2e^{(2L+C_{p,d,L}^2)t}(\epsilon_1 + \epsilon_2 + \epsilon_3C_p + I(N)), B = 2e^{(2L+C_{p,d,L}^2)t}L \text{ and } C=2e^{(2L+C_{p,d,L}^2)t}C_{p,d,L}\text{.}
\end{eqnarray*}

By using \cref{lemma3} once again,  we can obtain
\begin{eqnarray*}
\sup _{s  \in [0,t]}&\left\|W_p(\mu_s,\nu_s)\right \|_p  \leq 2e^{(2B+C^2)t} A\text{.}
\end{eqnarray*}

Taking this estimate back to $\varphi(t)$, and by (\ref{target}) we finally prove that there exist some constant $C_1>0$  such that
$$\left \| \sup _{u \in [0,t]}|X_u-Y_u| \right \|_p\leq C_1, $$
where $C_1$ depends on $L,T,N,d,\epsilon_1,\epsilon_2,\epsilon_3$ and will converges to $0$ as $N$ goes to infinity and  $\epsilon_1,\epsilon_2,\epsilon_3$ goes to $0$.

On the other hand, by the triangle inequality we have
\begin{eqnarray*}
\sup _{s  \in [0,t]}\left\|W_p(\mu_s,\hat\nu_s)\right \|_p  &\leq& \sup _{s  \in [0,t]}\left\|W_p(\mu_s,\nu_s)\right \|_p +\sup _{s  \in [0,t]}\left\|W_p(\nu_s,\hat\nu_s)\right \|_p \text{,}
\end{eqnarray*}
where the second term in the r.h.s of the inequality above  tends to $0$ as $N$ approaches to infinity due to  Lemma \ref{lemma4}, and hence the desired result follows.
\end{proof}

\subsection{Proof for Theorem 3.9}
\begin{proof}
By Using \cref{lemma1}, \cref{lemma2} and conditions, we have
\begin{small}
\begin{eqnarray}
&&\left\|\sup _{u \in [0,t]}|X_u-Y_u|\right \|_2 \leq
\left\| \sup _{u \in [0,t]}|X_0-Y_0
+\int_{0}^{u}(b(X_s,\hat\kappa_{s})-b ^{Pseudo}_{s})ds \right\|_2 \nonumber \\
&+&
\left\|\sup _{u \in [0,t]}\int_{0}^{u}(\sigma(X_s,\hat\kappa_{s})-\sigma ^{Pseudo}_{s})dB_s \right\|_2 \nonumber \\
&\leq &  \left\|\sup _{u \in [0,t]}\int_{0}^{u}(b(X_s,\hat\kappa_{s})- b(Y_s,\hat\varrho_{s}))ds \right\|_2\ +
\left\|\sup _{u \in [0,t]}\int_{0}^{u}(\sigma(X_s,\hat\kappa_{s})- \sigma(Y_s,\hat \varrho_{s}))ds \right\|_2\  \nonumber \\
&+& \epsilon_1 + \epsilon_2 + \epsilon_3C \nonumber \\
&\leq &  L\int_{0}^{t}(\left\| X_s- Y_s\right\|_2 + \left\| W_p(\hat\kappa_s,\hat\varrho_s)\right\|_2)ds  + C_{d}\sqrt{\int_{0}^{t}   \left \| \left \|\sigma(X_s,\hat \kappa_s ) - \sigma(Y_s,\hat \varrho_s ) \right\|  \right\|^2_2 ds} \nonumber \\
&+& \epsilon_1 + \epsilon_2 + \epsilon_3C \nonumber \\
&\leq & L\int_{0}^{t}(\left\| X_s- Y_s\right\|_2 + \left\| W_p(\hat\kappa_s,\hat\varrho_s)\right\|_2)ds + C_{d,L}\left \{\int_{0}^{t}(\left\| X_s- Y_s\right\|^2_2 + \left\| W_2(\hat\kappa_s,\hat\varrho_s)\right\|^2_2)ds\right \}^\frac{1}{2}  \nonumber \\
&+& \epsilon_1 + \epsilon_2 + \epsilon_3C. \nonumber
\end{eqnarray}
\end{small}
Define
\begin{eqnarray*}
 \varphi(t) := \epsilon_1 + \epsilon_2 + \epsilon_3C+
L\int_{0}^{t}\left\| W_2(\hat\kappa_s,\hat\varrho_s)\right\|_2ds
+ C_{d,L}\left \{\int_{0}^{t}\left\| W_2(\hat\kappa_s,\hat\varrho_s)\right\|^2_2ds\right\}^\frac{1}{2},
\end{eqnarray*}
then 
\begin{eqnarray}
\left\|\sup _{u \in [0,t]}|X_u-Y_u|\right \|_2  \leq
L\int_{0}^{t}\left\| X_u- Y_y\right\|_2ds
+C_{d,L}\left \{\int_{0}^{t}\left\|X_u- Y_u \right\|^2_2ds\right\}^\frac{1}{2} + \varphi(t).\nonumber
\end{eqnarray}
By \cref{lemma3}, we have
\begin{equation}\label{wq3}
\left \| \sup _{u \in [0,t]}|X_u-Y_u|
\right \|_2\leq 2e^{(2L+C_{d,L}^2)t}\varphi(t).
\end{equation}
On the other hand,
\begin{equation*}
\sup_{u  \in [0,t]}\mathbb{E} W_2^2(\hat\kappa_u,\hat\varrho_u)
\leq \mathbb{E} \left [\sup_{u\in [0,t]}|X_u-Y_u|^2\right ] = \left \| \sup _{u \in [0,t]}|X_u-Y_u|\right \|^2_2 
\leq (2e^{(2L+C_{d,L}^2)t}\varphi(t))^2.
\end{equation*}
Substituting the definition of $\varphi(t)$ in the inequality above, we get
\begin{eqnarray}
\sup_{s \in [0,t]}\left\|W_2(\hat\kappa_s,\hat\varrho_s)\right \|_2 &\leq&
2e^{(2L+C_{d,L}^2)t} \left ( \epsilon_1 + \epsilon_2 + \epsilon_3C+
L\int_{0}^{t}\left\|W_2^2(\hat\kappa_s,\hat\varrho_s)\right\|_pds
+ C_{d,L}\left \{\int_{0}^{t}\left\| W_2^2(\hat\kappa_s,\hat\varrho_s)\right\|^2_2ds\right\}^\frac{1}{2} \right )  \nonumber \\
&\leq& A+B\int_{0}^{t}\sup _{m \in [0,s]}\left\| W_2^2(\hat\kappa_s,\hat\varrho_s)\right\|_pds
+ C\left \{\int_{0}^{t}\sup _{m \in [0,s]}\left\|W_2^2(\hat\kappa_s,\hat\varrho_s)\right\|^2_2ds\right\}^\frac{1}{2} \nonumber
\end{eqnarray}
\begin{eqnarray*}
\text{where } A = 2e^{(2L+C_{d,L}^2)t}(\epsilon_1 + \epsilon_2 + \epsilon_3C) ,B = 2e^{(2L+C_{d,L}^2)t}L \text{ and } C=2e^{(2L+C_{d,L}^2)t}C_{d,L}\text{.}
\end{eqnarray*}
By using \cref{lemma3} we can obtain:
\begin{eqnarray}
\sup _{s  \in [0,t]}\left\|W_2(\hat\kappa_s,\hat\varrho_s)\right \|_2
\leq 2e^{(2B+C^2)t}A \text{.}
\end{eqnarray}
Taking the above estimate back to $\varphi(t)$ and then using \cref{wq3}, we obtain that there exist two constants $C_1$ and $C_2$ depending on $L,T,d,\epsilon_1,\epsilon_2,\epsilon_3$, such that
\begin{eqnarray}
\left \| \sup _{u \in [0,t]}|X_u-Y_u|
\right \|_2\leq
 C_1\text{,} \\
 \sup _{u  \in [0,t]}\left\|W_2(\hat\kappa_{u},\hat\varrho_{u})\right \|_2
\leq
 C_2\text{.}
\end{eqnarray}
And $C_1$ and $C_2$ will converge to $0$ when $T \to \infty$ and $\epsilon_1, \epsilon_2, \epsilon_3 \to 0 $.
\end{proof}

\end{document}

%% file: ex_shared.tex
\usepackage{lipsum}
\usepackage{amsfonts}
\usepackage{graphicx}
\usepackage{epstopdf}
\usepackage{algorithmic}
\ifpdf
  \DeclareGraphicsExtensions{.eps,.pdf,.png,.jpg}
\else
  \DeclareGraphicsExtensions{.eps}
\fi

\newsiamremark{remark}{Remark}
\newsiamremark{hypothesis}{Hypothesis}
\crefname{hypothesis}{Hypothesis}{Hypotheses}
\newsiamthm{claim}{Claim}

\headers{Solving McKean-Vlasov Equation by deep learning particle method}{Jingyuan Li, Wei Liu}

\title{Solving McKean-Vlasov Equation by deep learning particle method\thanks{Preprint, under review}}

\author{Jingyuan Li\thanks{School of Mathematics and Statistics, Wuhan University, Wuhan, Hubei PR China
  (\email{ljymath4@whu.edu.cn}}).
\and Wei Liu\thanks{School of Mathematics and Statistics, Wuhan University, Wuhan, Hubei PR China
  (\email{wliu.math@whu.edu.cn}}).}

\usepackage{amsopn}


%% file: ex_article.bbl
\begin{thebibliography}{10}

\bibitem{Antonelli2002RateOC}
{\sc F.~Antonelli and A.~Kohatsu}, {\em Rate of convergence of a particle
  method to the solution of the mckean-vlasov's equation}, Annals of Applied
  Probability, 12 (2002), pp.~423--476.

\bibitem{Baladron2012MeanfieldDA}
{\sc J.~Baladron, D.~Fasoli, O.~D. Faugeras, and J.~Touboul}, {\em Mean-field
  description and propagation of chaos in networks of hodgkin-huxley and
  fitzhugh-nagumo neurons}, Journal of Mathematical Neuroscience, 2 (2012),
  pp.~10 -- 10.

\bibitem{Bao2021FirstorderCO}
{\sc J.~Bao, C.~Reisinger, P.~Ren, and W.~Stockinger}, {\em First-order
  convergence of milstein schemes for mckean–vlasov equations and interacting
  particle systems}, Proceedings of the Royal Society A, 477 (2021).

\bibitem{Barbu2018FromNF}
{\sc V.~Barbu and M.~Rockner}, {\em From nonlinear fokker–planck equations to
  solutions of distribution dependent sde}, The Annals of Probability,  (2018).

\bibitem{Bogachev2019OnCT}
{\sc V.~I. Bogachev, M.~R{\"o}ckner, and S.~V. Shaposhnikov}, {\em On
  convergence to stationary distributions for solutions of nonlinear
  fokker–planck–kolmogorov equations}, Journal of Mathematical Sciences,
  242 (2019), pp.~69 -- 84.

\bibitem{Bossy1997ASP}
{\sc M.~Bossy and D.~Talay}, {\em A stochastic particle method for the
  mckean-vlasov and the burgers equation}, Math. Comput., 66 (1997),
  pp.~157--192.

\bibitem{carmona2018probabilistic}
{\sc R.~Carmona, F.~Delarue, et~al.}, {\em Probabilistic theory of mean field
  games with applications I-II}, Springer, 2018.

\bibitem{Chaintron2022PropagationOC}
{\sc L.-P. Chaintron and A.~Diez}, {\em Propagation of chaos: A review of
  models, methods and applications. i. models and methods}, Kinetic and Related
  Models,  (2022), \url{https://api.semanticscholar.org/CorpusID:247187507}.

\bibitem{du2023empirical}
{\sc K.~Du, Y.~Jiang, and J.~Li}, {\em Empirical approximation to invariant
  measures for mckean--vlasov processes: Mean-field interaction vs
  self-interaction}, Bernoulli, 29 (2023), pp.~2492--2518.

\bibitem{durmus2020elementary}
{\sc A.~Durmus, A.~Eberle, A.~Guillin, and R.~Zimmer}, {\em An elementary
  approach to uniform in time propagation of chaos}, Proceedings of the
  American Mathematical Society, 148 (2020), pp.~5387--5398.

\bibitem{Fournier2013OnTR}
{\sc N.~Fournier and A.~Guillin}, {\em On the rate of convergence in
  wasserstein distance of the empirical measure}, Probability Theory and
  Related Fields, 162 (2015), pp.~707--738.

\bibitem{guillin2021kinetic}
{\sc A.~Guillin, W.~Liu, L.~Wu, and C.~Zhang}, {\em The kinetic fokker-planck
  equation with mean field interaction}, Journal de Math{\'e}matiques Pures et
  Appliqu{\'e}es, 150 (2021), pp.~1--23.

\bibitem{guillin2022uniform}
{\sc A.~Guillin, W.~Liu, L.~Wu, and C.~Zhang}, {\em Uniform poincar{\'e} and
  logarithmic sobolev inequalities for mean field particle systems}, The Annals
  of Applied Probability, 32 (2022), pp.~1590--1614.

\bibitem{han2024learning}
{\sc J.~Han, R.~Hu, and J.~Long}, {\em Learning high-dimensional mckean--vlasov
  forward-backward stochastic differential equations with general distribution
  dependence}, SIAM Journal on Numerical Analysis, 62 (2024), pp.~1--24.

\bibitem{He2022AnEE}
{\sc J.~He, S.~Gao, W.~Zhan, and Q.~Guo}, {\em An explicit euler method for
  mckean-vlasov sdes driven by fractional brownian motion}, ArXiv,
  abs/2209.04574 (2022).

\bibitem{he2022explicit}
{\sc J.~He, S.~Gao, W.~Zhan, and Q.~Guo}, {\em An explicit euler method for
  mckean-vlasov sdes driven by fractional brownian motion}, arXiv preprint
  arXiv:2209.04574,  (2022).

\bibitem{Kloeden1977TheNS}
{\sc P.~E. Kloeden and E.~Platen}, {\em The numerical solution of stochastic
  differential equations}, The Journal of the Australian Mathematical Society.
  Series B. Applied Mathematics, 20 (1977), pp.~8 -- 12.

\bibitem{Lagaris1997ArtificialNN}
{\sc I.~E. Lagaris, A.~C. Likas, and D.~I. Fotiadis}, {\em Artificial neural
  networks for solving ordinary and partial differential equations}, IEEE
  transactions on neural networks, 9 5 (1997), pp.~987--1000.

\bibitem{lasry2018mean}
{\sc J.-M. Lasry and P.-L. Lions}, {\em Mean-field games with a major player},
  Comptes Rendus Mathematique, 356 (2018), pp.~886--890.

\bibitem{liu2021long}
{\sc W.~Liu, L.~Wu, and C.~Zhang}, {\em Long-time behaviors of mean-field
  interacting particle systems related to mckean--vlasov equations},
  Communications in Mathematical Physics, 387 (2021), pp.~179--214.

\bibitem{liu2019optimal}
{\sc Y.~Liu}, {\em Optimal Quantization: Limit Theorem, Clustering and
  Simulation of the McKean-Vlasov Equation}, PhD thesis, Sorbonne
  universit{\'e}, 2019.

\bibitem{Liu2022ParticleMA}
{\sc Y.~Liu}, {\em Particle method and quantization-based schemes for the
  simulation of the mckean-vlasov equation}, ArXiv, abs/2212.14853 (2022).

\bibitem{Muniandy2001ModelingOL}
{\sc S.~V. Muniandy and S.~C. Lim}, {\em Modeling of locally self-similar
  processes using multifractional brownian motion of riemann-liouville type.},
  Physical review. E, Statistical, nonlinear, and soft matter physics, 63 4 Pt
  2 (2001), p.~046104.

\bibitem{Rached2022MultilevelIS}
{\sc N.~B. Rached, A.-L. Haji-Ali, S.~Mohan, and R.~Tempone}, {\em Multilevel
  importance sampling for mckean-vlasov stochastic differential equation},
  ArXiv, abs/2208.03225 (2022).

\bibitem{rached2022single}
{\sc N.~B. Rached, A.-L. Haji-Ali, S.~M.~S. Pillai, and R.~Tempone}, {\em
  Single level importance sampling for mckean-vlasov stochastic differential
  equation}, arXiv preprint arXiv:2207.06926,  (2022).

\bibitem{Raissi2019PhysicsinformedNN}
{\sc M.~Raissi, P.~Perdikaris, and G.~E. Karniadakis}, {\em Physics-informed
  neural networks: A deep learning framework for solving forward and inverse
  problems involving nonlinear partial differential equations}, J. Comput.
  Phys., 378 (2019), pp.~686--707.

\bibitem{reisinger2023convergence}
{\sc C.~Reisinger and M.~O. Tsianni}, {\em Convergence of the euler--maruyama
  particle scheme for a regularised mckean--vlasov equation arising from the
  calibration of local-stochastic volatility models}, arXiv preprint
  arXiv:2302.00434,  (2023).

\bibitem{Ryck2021OnTA}
{\sc T.~D. Ryck, S.~Lanthaler, and S.~Mishra}, {\em On the approximation of
  functions by tanh neural networks}, Neural networks : the official journal of
  the International Neural Network Society, 143 (2021), pp.~732--750.

\bibitem{Salinas2020ATF}
{\sc D.~G. Salinas and M.~O. Gallardo}, {\em A topic for integrated teaching of
  mathematics and biology: the parabola of chaos in tumour cell aneuploidy},
  International Journal of Mathematical Education in Science and Technology, 53
  (2020), pp.~909 -- 919.

\bibitem{sznitman1991topics}
{\sc A.-S. Sznitman}, {\em Topics in propagation of chaos}, Ecole
  d’{\'e}t{\'e} de probabilit{\'e}s de Saint-Flour XIX—1989, 1464 (1991),
  pp.~165--251.

\bibitem{veretennikov2006ergodic}
{\sc A.~Y. Veretennikov}, {\em On ergodic measures for mckean-vlasov stochastic
  equations}, in Monte Carlo and Quasi-Monte Carlo Methods 2004, Springer,
  2006, pp.~471--486.

\bibitem{wang2018distribution}
{\sc F.-Y. Wang}, {\em Distribution dependent sdes for landau type equations},
  Stochastic Processes and their Applications, 128 (2018), pp.~595--621.

\end{thebibliography}
